\documentclass[a4paper]{article}
\usepackage{mathtools,amsmath,amssymb,amsthm,amsfonts,graphicx,fullpage}

\usepackage{changepage}

\usepackage[utf8x]{inputenc}

\usepackage{textcomp,marvosym}

\usepackage{cite}

\usepackage{nameref,hyperref}

\usepackage[table]{xcolor}

\usepackage{array}

\newcolumntype{+}{!{\vrule width 2pt}}

\usepackage[outdir=./]{epstopdf}

\newcommand{\comment}[1]{}

\newcommand{\mR}{\mathcal{R}}
\newcommand{\N}{\mathcal{N}}
\newcommand{\rewieverTwo}{\textcolor{black}}
\newcommand{\rewieverOne}{\textcolor{black}}
\newcommand{\yves}{\textcolor{black}}

\newtheorem{remark}{Remark}

\usepackage{lscape}

\makeatother

\title{Modeling the impact of rainfall and temperature on sterile insect control strategies in a Tropical environment} 
\author{Dumont Yves$^{1,2,3}$\footnote{Corresponding author: yves.dumont@cirad.fr}, Duprez Michel$^{4}$ \\
\small $^{1}$ CIRAD, UMR AMAP, 97410 St Pierre, R\'eunion Island - France \\
\small $^{2}$ AMAP, Univ Montpellier, CIRAD, CNRS, INRAE, IRD, Montpellier, France \\
\small $^{3}$ Department of Mathematics and Applied Mathematics, University of Pretoria, Pretoria, South Africa\\
\small $^{4}$ Inria, \'equipe  Mimesis,  Universit\'e  de  Strasbourg,  Icube, 1  Place  de  l'H\^opital,  Strasbourg,  France
}

\begin{document}

\maketitle
\begin{abstract}
The sterile insect technique (SIT) is a biological control technique that can be used either to eliminate or decay a wild mosquito population under a given threshold to reduce the nuisance or the epidemiological risk. In this work, we propose a model using a differential system that takes into account the variations of rainfall and temperature over time and study their impacts on sterile males' releases strategies. Our model is as simple as possible to avoid complexity while being able to capture the temporal variations of an \textit{Aedes albopictus} mosquito population in a domain treated by SIT, located in R\'eunion island. The main objective is to determine what period of the year is the most suitable to start a SIT control to minimize the duration of massive releases and the number of sterile males to release, either to reduce the mosquito nuisance, or to reduce the epidemiological risk. Since sterilization is not $100\%$ efficient, we also study the impact of different levels of residual fertility within the released sterile males population. Our study shows that rainfall plays a major role in the dynamics of the mosquito and the SIT control, that the best period to start a massive SIT treatment lasts from July to December, that residual fertility has to be as small as possible, at least for nuisance reduction. \yves{Indeed, when the main objective is to reduce the epidemiological risk, we show that residual fertility is not necessarily an issue.} Increasing the size of the releases is not always interesting. We also highlight the importance of combining SIT with mechanical control, i.e. the removal of breeding sites, in particular when the initial mosquito population is large. \yves{Last but not least our study shows the usefulness of the modeling approach to derive various simulations to anticipate issues and demand in terms of sterile insects' production.}
\end{abstract}



\section{Introduction}
Being vectors of many diseases, like malaria, dengue, Lymphatic filariasis, zika, chikungunya, yellow  fever, and japanese  encephalitis, mosquitoes are one of the biggest killers in the world, and particularly in tropical and subtropical areas. Indeed, Female mosquitoes need blood meals to provide the nutriments for egg development. That is why they bite either during the night, like \textit{anopheles spp}, at sunset or sunrise, like \textit{aedes albopictus}, or along the day, like \textit{aedes aegypti}, inside or outside houses. Since they can survive several weeks, they will bite several times depending on their gonotrophic cycle \cite{Klowden1994}. While opportunistic feeders, Anopheles and Aedes mosquitoes prefer biting mammals and, preferably, humans \cite{Tirados2006}, with sometimes multiple blood feeding \cite{Delatte2010}.


To eradicate mosquito-borne diseases, the initial option was the massive use of chemicals in the 1950's and 1960's, mainly against anopheles. Despite some successes, like in R\'eunion island where malaria has been eradicated since the fifties, thanks to the use of DDT \cite{Hamon1954}, we know now that this was a huge mistake: even if the use of adulticides was successful against mosquitoes, the damages on the biodiversity  were important. In addition, resistance to some adulticides raised, such that, in several places around the World, there is no option left to fight the mosquitoes, like in the French West Indies \cite{Marcombe2012}. Fortunately, other eco-friendly control methods, more respectful of the biodiversity, have been developed. Among them, the Sterile Insect Technique (SIT) is the most promising one. This is an old control technique, proposed in the 30s and 40s by three key researchers in the USSR, Tanzania and the USA and, first, applied in the field in the 50's \cite{SIT}. SIT has been used more or less successfully on the field against various kinds of Pests or Vectors (see \cite{SIT,Aronna2020} for various examples). The classical SIT consists of mass releases of males sterilized by ionizing radiation. The released sterile males transfer their sterile sperms to wild females, which will have no viable offspring, resulting in a progressive decay of the targeted population. For mosquitoes, other sterilization techniques have been developed using either genetics (release of insects carrying a dominant lethal technique, in short,  RIDL technique \cite{Solari2020}) or bacteria (SIT-IIT : cytoplasmic incompatibility thanks to \textit{Wolbachia} bacteria) \cite{Sinkins2004}. Many SIT projects against mosquitoes are now ongoing around the World. \rewieverOne{It is important to keep in mind that controlling a vector population could have an impact on the ecosystem equilibrium. In fact, it depends if the targeted vector is an invasive species, i.e. established recently or if it has been established for a long time. If the establishment is recent (a couple of years), then with elimination, we can expect to recover the initial state of the ecosystem before the invasion starts. On the contrary, when the vector is established, then elimination is not necessarily the goal. Reduction is preferable in order to reduce the epidemiological risk. We will explore these two objectives in this work.}

Several mosquito models have been developed, from very simple models to more sophisticated ones, according to the number of biological states that are taken into account. These models consider the mosquito population either at the individual level (see for instance \cite{Almeida2010}) or at the population level (see for instance \cite{Li2011}), taking eventually into account the spatial component (see for instance \cite{Anguelov2020_arxiv} for a brief review). Each of these approaches has advantages and disadvantages. However, it is always important to keep in mind that biological parameters are difficult to obtain in the field, such that for too complex models, it might be impossible to set the (right) values of the parameters. In addition, most of the biological parameters are linked to environmental parameters, like rainfall, humidity and temperature. 
The main goal of the present article will be to take into account these parameters in the mathematical model to find the best control.

Mathematically, it can make sense to consider constant or periodic parameters to derive the main dynamics of the mosquito population. However, for practical reasons, i.e. anticipating when the mosquito population is growing or decaying, it is better to have a mosquito model that takes into account the most important parameters, like rainfall, humidity and temperature, \yves{to derive simulations }to adapt the control, i.e. for instance the size of the releases or/and the period of the releases. \yves{Las but not least, being able to derive several realistic simulations can help to anticipate the demand in terms of sterile insects production.}

Since our work takes place within an ongoing SIT feasibility program (TIS 2B) against \textit{Aedes albopictus} in La R\'eunion, a tropical French island located in the Indian Ocean, we will focus on this species in the rest of the paper. In La R\'eunion, \textit{Aedes albopictus} has become the main vector of dengue and chikungunya \cite{Vazeille2018}. However, we believe that our approach is sufficiently generic to be applied to other mosquito species. 

While links between \textit{Aedes albopictus } parameters and temperature have been studied in laboratory \cite{Delatte2009}, for a fixed (laboratory) humidity, it is more challenging to take into account the effect of rainfall (and humidity).

Some weather-dependent mosquito models have been developed, mainly with Temperature-dependent parameters (see for instance \cite{Dufourd2013,Anopheles2015} and references therein) and very few with temperature and rainfall-dependent parameters (see \cite{Tran2020} and references therein). However, in general, these last models are quite complex: 
they relied on statistical approaches and on the user's subjective choices, such that the calibration (of many parameters), with respect to the environmental parameters, is not generic and might not be able to provide a unique set of valuable values. We firmly believe that simple (but not too simple) models can rapidly provide useful and reliable information to help field experts to manage vector control campaigns.

That is why we consider a minimalistic model (minimal in terms of stages and thus parameters), based on \cite{AnguelovTIS2020}, to build an approximation of the mosquito population over many years. Then, we show that the density variation (in time) of the mosquito population can help to improve vector-control strategies combining SIT with mechanical control, i.e. the removal of breeding sites to reduce the larvae carrying capacity. 

It is well known that sterilization does not necessarily induce that the sperm of the sterilized males is $100\%$ sterile (see for instance \cite{oliva2012sterile,Iyaloo2020}): it depends on the radiation source, the dose-rate, and the container. That is why quality control after sterilization is set up to control that the sterile males remain competitive and efficient. Efficiency is thus related to the effective sterility or the residual fertility (RF, shortly) , i.e. the percentage of sperms that is still fertile despite the irradiation, that is equivalent to say, at the population level, that a small proportion, $\varepsilon$, of the sterile males is fertile.
In our model, we will take into account  the fact that residual fertility can occur within the sterile 
male population to study its impact on SIT control. 

In general, in most of the SIT models (except in \cite{Aronna2020} and the present paper), residual fertility is never taken into account, while it can have a strong impact on SIT efficiency.
Indeed, in \cite{Aronna2020}, using a 2-dimensional mosquito model with constant parameters, the authors showed that the residual fertility, $\varepsilon$, is strongly linked to the basic offspring number, $\mathcal{N}$, of the wild population: when
\begin{equation}
\label{condeps}
\varepsilon<\dfrac{1}{\N},
\end{equation}
then SIT is efficient and the mosquito population can be lowered under any given threshold, provided that enough sterile males are released. However, when $\varepsilon$ is below but close to $\dfrac{1}{\N}$, the amount of sterile males increases almost exponentially \cite{Aronna2020}.
On the contrary, when $\varepsilon>\dfrac{1}{\N}$, then, whatever the size of the sterile males 
releases, the wild population cannot be lowered under a certain threshold value, that can be roughly estimated \cite{Aronna2020}. In fact, we suspect that a too-large residual fertility could partly explain failures in some SIT programs. \rewieverTwo{For instance, in \cite{Iyaloo2020}, a SIT program conducted on Mauritius island, which is close to R\'eunion island, the residual fertility was experimentally estimated at around $3.05\%$ $(2.29\%-3.92\%)$, which could explain why mitigated results were obtained, before a cyclone occurred and broke the experiment. We believe that $3\%$ of RF is too large such that the reduction of the wild population with SIT is limited. Indeed, in \cite{Aronna2020}, the authors showed that if RF is too large, then, whatever the size of the sterile males releases, the population cannot be lowered under a given threshold.}

The paper is organized as follows: In section 2, we build a temperature and rainfall-dependent entomological model and an impulsive periodic SIT model. Then, in section 3, we provide several numerical simulations to discuss the impact of the temperature, the rainfall, the residual fertility and the mechanical control on SIT starting period and duration to reduce either the nuisance or the epidemiological. Finally, in section 4, we end the paper with some conclusions and perspectives.

\section{Rainfall, humidity, and temperature dependent SIT model}
The first aim of the present work is to develop a temperature-rainfall-dependent entomological model, to take into account real field data, including mean daily rainfall and mean daily temperature.

Following \cite{AnguelovTIS2020}, we will consider the following model
\begin{equation}\label{eq:no SIT}
\left\{ \begin{array}{l}
\dfrac{dA}{dt}=\phi(Temp)F-\left(\gamma(Temp)+\mu_{1,A}(Temp)+\mu_{A,2}(Temp,Rain)A\right)A,\\ \\
\dfrac{dM}{dt}=(1-r(Temp))\gamma(Temp)A-\mu_{M}(Temp)M,\\\\
\dfrac{dF}{dt}=r(Temp)\gamma(Temp)A-\mu_{F}(Temp)F,
\end{array}\right.
\end{equation}
where $A$, $M$, and $F$ represent, respectively the aquatic (larvae, pupae) stage, the adult (male and \rewieverTwo{mated and fertile }female) stages. \rewieverTwo{Here, we implicitly assume that all emerging females will mate and become fertile. Note carefully that compartment $A$ does not include the eggs since we consider a density mortality rate $\mu_{A,2}$ that concerned only the larval and pupal stages. Thus, the parameter $\phi$ represents the daily amount of deposited eggs that hatch to become larvae. All parameters of model \eqref{eq:no SIT} are described in Table \ref{table:model2}.}
\begin{table}[h!]
\begin{center} 
  \caption{\bf Description of the entomological parameters \label{table:model2}}
 \begin{tabular}{|l|l|l|}
    \hline
    Symbol & Description & Unit\\
    \hline \hline \hline
    $\phi$ & Number of \rewieverTwo{hatched} eggs at each deposit per capita   & Day$^{-1}$ \\
    \hline
    $\gamma$ & Maturation rate from larvae to adult  & Day$^{-1}$ \\
    \hline
    $\mu_{A,1}$ & Density independent mortality rate of the aquatic stage  & Day$^{-1}$ \\
    \hline
    $\mu_{A,2}$ & Density dependent mortality rate of the aquatic stage  & Day$^{-1}$ Individuals $^{-1}$ \\
    \hline
    $r$ & Sex ratio & - \\
     \hline
     $1/\nu_m$ & Average extrinsic incubation period (EIP) & Day \\
     \hline
    $1/\mu_F$ & Average lifespan of fertilized and eggs-laying females & Day \\
         \hline
    $1/\mu_M$ & Average lifespan of males & Day \\
         \hline
    \hline
  \end{tabular}
  \end{center}
\end{table}
\rewieverTwo{This model has been studied in \cite{AnguelovTIS2020} when the parameters are supposed to be constant. In \cite{AnguelovTIS2020}}the authors developed a new strategy to maintain the wild population under a certain threshold, using a permanent and sustainable low level of SIT control, thanks to a  massive-small releases strategy. We set $\N_{\max}=\max_{t\in [0,+\infty)}\N(t)$, where $\N$ represents the basic offspring number, defined as follows
$$
\N=\dfrac{r\phi\gamma}{\left(\gamma+\mu_{A,1}\right)\mu_F}.
$$

In R\'eunion island, a network of Weather stations (from M\'et\'eo France but also from CIRAD) allows us to estimate some weather parameters all around the island, and in particular where real SIT field experiments started in July 2021, in the site of Duparc, a 20-hectare urban area located within the commune of Sainte Marie in the northern district of La R\'eunion \cite{Legoff2019}. These releases consist of manually weekly releases of 150 000 to 250 000 sterile males (once a week) produced and irradiated in La R\'eunion. The efficacy of sterile males releases is assessed over time by monitoring the induced sterility in eggs using ovitraps and the subsequent population suppression using adult traps, and, from time to time, Mark-Release-Recapture experiments. This experiment lasted until September 2022.

From the Weather station located in La Mare, close to the site of Duparc (1km), we are able to obtain the following daily data: the rainfall, the average daily temperature and also the humidity. Thus, following \cite{Valdez2018}, we first define the breeding site carrying capacity to define the density death-rate $\mu_{A,2}$. Obviously, the persistence of breeding sites is a key factor for the mosquito population's survival. Indeed, rainfall creates breeding sites, while evapotranspiration tends to shrink them. Following \cite{Valdez2018}, we define the variable $H(t)$ as the amount of water available at
day $t$ and defined as follows
\[
H(t+1)=\left\{ \begin{array}{ll}
0 & \text{if }H(t)+\Delta(t)\leq0,\\
H_{max} & \text{if }H(t)+\Delta(t)\geq H_{max},\\
H(t)+\Delta(t) & \text{otherwise},
\end{array}\right.
\]
with 
\[H_{\mbox{\footnotesize max}}:=\max_{t\in [t_0,t_0+T]}\{Rain(t)\}
\]
and
\[
\Delta(t):=Rain(t)-Evap(t),
\]
where $Rain(t)$ is the daily rainfall and $Evap(t)$ the
daily evaporation, $t_0$ is the starting date of the simulations, and $T$ the total duration. Following \cite{Valdez2018}, the evaporation function is
defined as follows
\[
Evap(t)=k\times\left(25+Temp(t)^{2}\right)\times\left(100-Hum(t)\right),
\]
where $Temp(t)$ is the average temperature and $Hum(t)$ the humidity.
Finally, the carrying capacity is defined 
\[
K(t)=K_{max}\dfrac{H(t)}{H_{max}}+K_{0},
\]
where $K_{0}>0$ can be seen as the fixed artificial carrying capacity, i.e.
rainfall-independent, human-made, by watering, for instance, plants around houses (flower
pots, plates, and vases), and $K_{\max}$, the natural maximal carrying capacity.

The initialization in time of the variables (in particular of $H$) is explained in the simulations section, page \pageref{sec:simu}. We notice that, in general, in the literature, the question of initialization is not always taken into account, while it can have a strong influence over several months on the population dynamic. We illustrate this fact in Figure \ref{fig:init H}, page \pageref{fig:init H}. 

In order to estimate $\mu_{2,A}$, we consider the positive equilibrium related to the carrying capacity, like in \cite{Dumont2012,Dufourd2013,Strugarek2019}. Thus, for a fixed value of $K$, the aquatic stage at equilibrium is given by
\begin{equation}
A_{K}^{*}=\left(1-\dfrac{1}{\mathcal{N}}\right)K.
\label{eq1}    
\end{equation}

In our model, the aquatic stage at equilibrium is defined as follows
\begin{equation}
    \label{eq2}
A^{*}=\dfrac{\gamma+\mu_{A,1}}{\mu_{A,2}}\left(\mathcal{N}-1\right),
\end{equation}
such that considering the equality between both equilibria given in \eqref{eq1} and \eqref{eq2}, we derive the following relationship between $\mu_{A,2}$ and $K$, that
is
\begin{equation*}
\mu_{A,2}=\dfrac{\mathcal{N}(\gamma+\mu_{A,1})}{K}=\dfrac{r\gamma\phi}{\mu_{F}K}.
\end{equation*}
Thus, taking into account the dependency of the parameters to Temperature and Rainfall, we deduce that 
\begin{equation}
\mu_{A,2}(Temp,Rain)=\dfrac{r(Temp)\gamma(Temp)\phi(Temp)}{\mu_{F}(Temp)K(Temp,Rain)}.
   \label{muA2}
\end{equation}
As it is well known, eggs are deposited above the waterline. They
hatch once they are flooded by rainfall. However, the hatching rate
can be seasonal. From \textit{A. albopictus} eggs, picked up at Duparc, the mean hatching proportion is around $90\%$ (G. Legoff, personal communication, TIS2B project).

For the other parameters, we will consider the data recalled \yves{in Appendix A, page \pageref{AppendixA},} in Table \ref{tab:parameters1}, page \pageref{tab:parameters1}, as used in \cite{Strugarek2019} to obtain the 
parameters estimates given in Table \ref{tab:parameters2}, page \pageref{tab:parameters2}.

Then, we can use simple interpolation polynomials (with  cubic spline, like in \cite{Dufourd2012B}) to estimate these parameters for any given temperature
$Temp$, within the range $[15^{o},35^{o}]$. 

We will assume constant or periodic releases of sterile males, at rate $u_S(t)$ (either constant or variable), and $\mu_S$, the sterile male mortality rate, which  is supposed to be similar to the wild males mortality rate \cite{Oliva2012}.
The dynamic of the sterile males is modeled by
\begin{equation}\label{eq:MS}
\dfrac{dM_{S}}{dt}=u_S(t)-\mu_{S}(Temp)M_{S}.
\end{equation}

\comment{
Then following \cite{Aronna2020}, we will consider the residual fertility, $\varepsilon$, in our SIT model. As recalled in the introduction, for constant value parameters, we need to verify  $\varepsilon<\dfrac{1}{\N}$ in order to have an efficient SIT control, i.e. such that the wild population can be lowered under a given value, using appropriate (massive) releases. In \cite{Aronna2020}, using a minimalistic model with constant parameters, the authors showed that, for instance, for \textit{A. albopictus} mosquito, the RF should be lower than $2.5\%$. 
}

\rewieverTwo{Thus, based on model \eqref{eq:no SIT}, we consider} the following temperature and rainfall dependent SIT-model
\begin{equation}
\left\{ \begin{array}{rcl}
\dfrac{dA}{dt}&=&\phi(Temp)F-\left(\gamma(Temp)+\mu_{1,A}(Temp)+\mu_{A,2}(Temp,Rain)A\right)A,\smallskip\\\smallskip
\dfrac{dM}{dt}&=&(1-r(Temp))\gamma(Temp)A-\mu_{M}(Temp)M,\\
\dfrac{dF}{dt}&=&r(Temp)\gamma(Temp)\dfrac{M+\varepsilon \beta M_S}{M+\beta M_{S}}A-\mu_{F}(Temp)F, \\
\dfrac{dM_{S}}{dt}&=&u_S(t)-\mu_{S}(Temp)M_{S}.
\end{array}\right.
\label{TIS_mod}
\end{equation}
\rewieverTwo{The parameter $\beta$ is the competition parameter that is deduced from the Fried index \cite{Fried}. Following \cite{oliva2012sterile}, we will consider $\beta=1$, which means that a sterile male is as competitive as a wild male. The Fried index is estimated through a protocol described in \cite{Fried} (see also \cite{WHO_IAEA_2020}). Last, as explained in the introduction, all irradiated males are not $100\%$ sterile: a small proportion of their sperm can stay fertile. At the population scale, this is modeled by the parameter $\varepsilon$, which represents the proportion of sterile males that remain fertile and thus can fertilize emerging females. As already explained in \cite{DUMONT2022}[Remark 4], \textit{A. albopictus} females got $3$ spermathecae that allow to stock sperms such that, in general, one mating is sufficient to fulfill at least $2$ spermathecae and to have enough sperms to fertilize eggs along its lifespan. Thus the nonlinear term related to SIT treatment occurs when females emerge from compartment $A$ and not in the birth rate. Considering SIT impact on the birth rate, i.e. $\dfrac{M+\varepsilon \beta M_S}{M+\beta M_S}$, is very convenient from the mathematical point of view, but it is not biologically realistic. We prefer to stick to the biological reality, at least for \textit{A. albopictus}. Thus, altogether, in model \eqref{TIS_mod}, the term $\dfrac{M+\varepsilon \beta M_S}{M+\beta M_S}$ represents the probability for an emerging female to be fertilized and thus to enter compartment $F$. 
}

\subsection{Impulsive SIT massive-small releases strategy }

\rewieverTwo{In the rest of this paper, we will consider the so-called massive-small releases strategy developed in \cite{AnguelovTIS2020}. Indeed, once a few sterile males are released, this induces a strong Allee effect, such that the elimination equilibrium, $\mathbf{0}$, and a positive equilibrium $E$ are both asymptotically stable. In \cite{AnguelovTIS2020}, we showed that it is possible to estimate the minimal amount of sterile males to release to stay in the basin of attraction of $\mathbf{0}$ when the wild population is sufficiently small. This means that massive releases have first to be used over a sufficient time period in order to enter in a particular subset of the basin of attraction of $\mathbf{0}$, that is estimated thanks to a given (small) release rate value. Once inside the basin of attraction of $\mathbf{0}$, it is possible to switch from massive releases to small releases. This strategy is particularly useful since SIT cannot be stopped in order to prevent the re-establishment of the targeted vector, either by introduction or by migration.}

\noindent \rewieverTwo{According to this massive-small releases strategy,} we want to find the best period to start SIT control (without or with  Mechanical Control) to minimize the number of massive periodic impulsive releases, and thus the amount of sterile males to release.

\noindent \yves{We provide some computations of the SIT-equilibria in Appendix B, page \pageref{AppendixB}.} Thanks to \cite{AnguelovTIS2020}, \rewieverOne{and assuming that $\N \varepsilon \leq 1$, }system  \eqref{TIS_mod} has the following long term behavior
\begin{itemize}
    \item There exists a release rate threshold, $u^*$, such that when $u_S(t)>u^*$ then $(A,M,F)$ converges to $\mathbf{0}$.
    \item When $0<u_S(t)<u^*$, then there exist two positive equilibria, $\bf{E_1}$ and $\bf{E_2}$ such that $\bf{E_1}<\bf{E_2}$. Moreover $\bf{E_1}$ is unstable while $\bf{E_2}$ is stable and $[\bf{0},\bf{E_1})$ lies in the basin of attraction of $\bf{0}$.
\end{itemize}

\begin{remark}
    \rewieverTwo{When $\N \varepsilon > 1$, then, elimination is not possible, whatever the size of the massive releases. Using the equilibrium value for the aquatic stage, $A$, given in appendix B, page \pageref{AppendixB}, it is straightforward to show that there exists a lower bound for the SIT positive equilibrium, for all $u_S>0$, given by 
    \begin{equation} 
    \left\{\begin{array}{l}
    A_l^*=2\dfrac{\varepsilon\mathcal{N}-1}{\mathcal{N}-1} A^*,\\
    M_l^*=\dfrac{(1-r)\gamma}{\mu_M}A_l^*,\\
    F_l^*=\dfrac{\gamma+\mu_{1,A}+\mu_{A,2}A_l^*}{\phi}A_l^*.
    \end{array}\right.
    \label{lowerbound}
    \end{equation}
    However, even if $\mathcal{N}\varepsilon>1$ and nuisance reduction is not possible, it does not mean that it is not possible to reduce the epidemiological risk. The lower bound \eqref{lowerbound} will be useful in the epidemiological risk section.}
\end{remark}

\comment{
\begin{remark}
    \rewieverOne{The SIT threshold parameter, $u^*$, is defined as follows
    \begin{equation}\label{seuil-ms}
			\begin{array}{lcl}
				u^{*} = \dfrac{\mu_{M_S}}{\beta \mathcal{Q}}\left(\sqrt{\mathcal{N}\left(1-\varepsilon\right)}-\sqrt{1-\mathcal{N}\varepsilon}\right)^{2},
			\end{array}
		\end{equation}
    where $\mathcal{Q}=\dfrac{\mu_{A,2}\mu_{M}}{(\gamma+\mu_{A,1})(1-r)\gamma}$.}
\end{remark}
}
We will consider two different levels of massive releases  ($6000$ or $12000$ sterile males per ha), \rewieverOne{such that $u>u^*$ for all $t\geq 0$,} to decrease the wild population below a threshold set, $[\bf{0},\bf{E_1})$, defined by (inexpensive) small releases (say 100 sterile males per ha), in a minimum time.
We want to know how the duration of massive releases is influenced by temperate and rainfall conditions during the year, the residual fertility and the level of mechanical control.
In other words, we search a control with the form
\begin{equation}
    u(t)=\tau\Lambda_{\mbox{\footnotesize massive}}\sum_{i=1}^{N}\delta_{t_0+(i-1)\tau}(t),
    \label{u}
\end{equation}

where $\delta(t)$ is the Dirac function, $t_0$ is the starting time of the massive releases, $\tau$ the periodicity of the releases (here, $\tau=7$), and $N$ the number of weakly massive releases.
The massive release $\Lambda_{\mbox{massive}}$ and the small releases $\Lambda_{\mbox{\footnotesize small}}$ are fixed. For a given $t_0$ there exists a $t_1(t_0)=t_0+N_1(t_0)\tau$ (according to \cite{AnguelovTIS2020}) such that, after this time, the  wild population $(A,M,F)$ remains in the box
$[\bf{0},\bf{E_{1,min}}(\tau\Lambda_{small}))$ where $\bf{E_{1,min}}(\tau\Lambda_{small})$ is defined as follow:
For a given $\Lambda_{\mbox{\footnotesize small}}$, for each time $t$, we compute the equilibrium $\bf{E_1}(t)$ of the system associated to the parameters at time $t$, and we define  $\bf{E_{1,min}}(\tau\Lambda_{small})=\min_t \bf{E_1}(t),$
where the minimum is taken between the beginning and the end of the time interval considered in the simulation.

The main goal is to find the (best) starting time of the massive releases, $t_0$, in order to minimize the duration of the releases and thus, the number of massive releases.

System \eqref{eq:MS}-\eqref{u} can be rewritten as an impulsive differential system with fixed moments of impulse effect, that is \begin{equation}
    \left\{
    \begin{array}{l}
         \dfrac{dM_S}{dt}=-\mu_S(T)M_S,\qquad t_0+i\tau < t\leq t_0+(i+1)\tau  \\
         M(t^+)=M(t)+\tau\Lambda_{massive}, \qquad t=t_0+i\tau, 
    \end{array}
    \right.
    \label{impulsive}
\end{equation}
for $i=0,...,N-1$. Since, the right-hand side of \eqref{TIS_mod}-\eqref{impulsive}$_1$ is locally Lipschitz-continuous on $\mathbb{R}^4$, we can use a classical existence Theorem (for instance Theorem 1.1 in \cite{Bainov1995}, or Theorem 2.1 in \cite{Bainov1993}), to deduce that there exists $T_e>0$ and a unique solution of system \eqref{TIS_mod}-\eqref{impulsive}, defined from $(t_0,T_e)\longrightarrow \mathbb{R}^4$.

\section{Numerical simulations and discussion}\label{sec:simu}

System \eqref{TIS_mod}-\eqref{impulsive} is solved thanks to \texttt{odeint} of the python library \texttt{scipy.integrate}. The codes and the data corresponding to the temperature and the rainfall in R\'eunion island are available on 
\begin{multline*}
    \texttt{https://github.com/michelduprez/Impact-of}\\\texttt{-Rainfall-and-Temperature-on-IT-control-strategies.git}
\end{multline*}
\yves{Computations are not complex but very long. That is why we ran our code on the MESO@LR-Platform (University of Montpellier). The figures have been draw using the software Matlab \cite{MATLAB2019}.}

We consider temperature, rainfall and humidity data based on (noised) data recorded in weather stations close to the site of Duparc (a 20 ha place), a neighborhood of Sainte Marie, located in the North of R\'eunion island from the 1st of January 2009 to the 14th of July 2021.

In La R\'eunion, since the sterile males are produced on-site, everything from the eggs to the sterilization is controlled (qualitatively), such that since the production started, the residual fertility is, in general, less than $1\%$, with an average value of around $0.6\%$, a very good results compared to other SIT projects, like \cite{Iyaloo2020}. That is why, in the forthcoming simulations, we will consider three cases of RF, namely $0\%$, $0.6\%$ and $1.2 \%$.

\begin{figure}[!ht]
\centering
\begin{center}
\includegraphics[width=1.0 \linewidth]{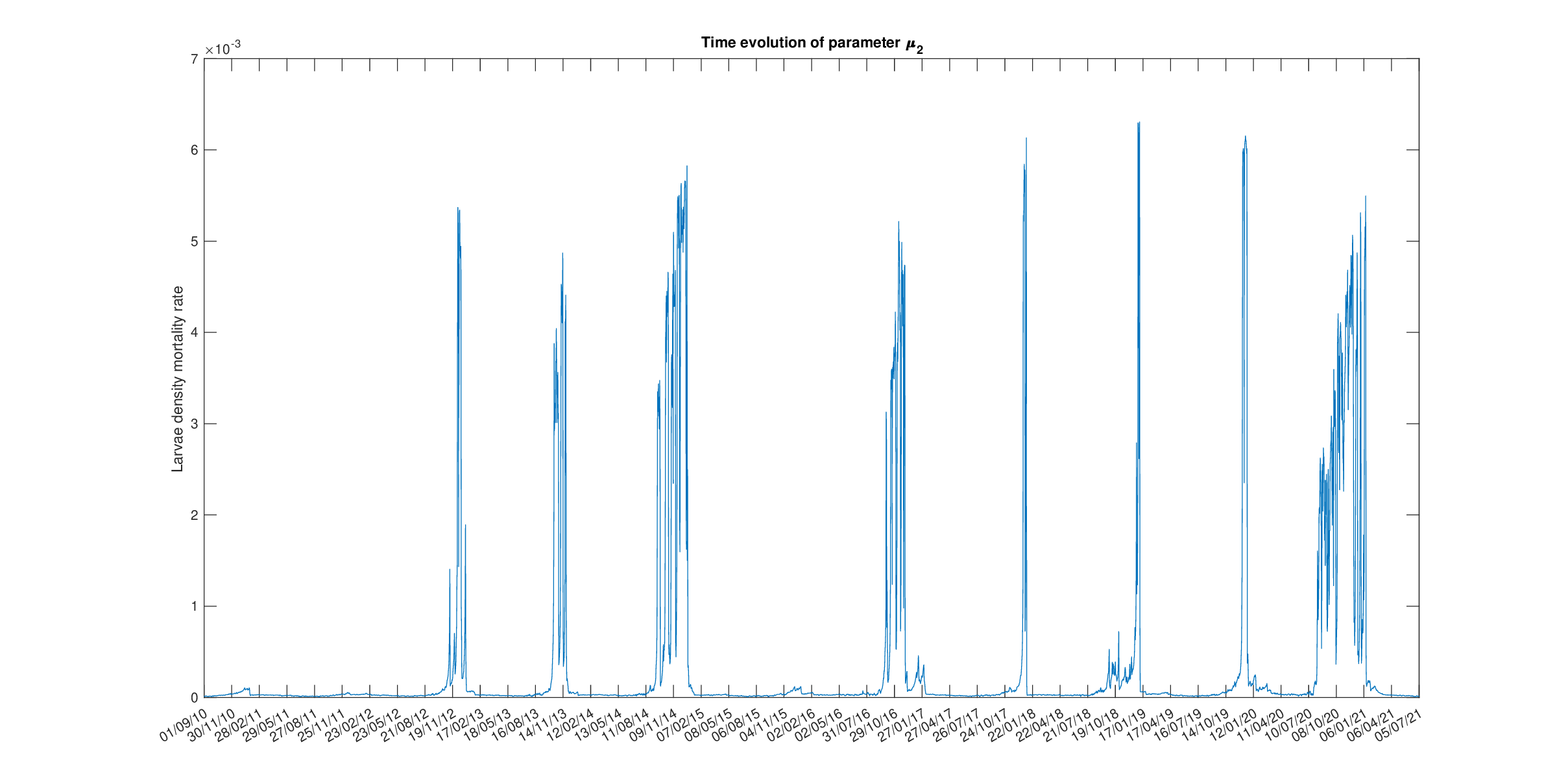} 
\end{center}
\caption{Time evolution of $\mu_{A,2}$ from the 1st of September 2010 until the 14th of July 2021}  \label{fig:1}
\end{figure}

In Fig. \ref{fig:1}, we consider the following values: $K_{\max}=20\times 10000$,  
$K_0=20 \times 100$ to derive the carrying capacity, $K$, and thus $\mu_{A,2}$.

It is interesting to see the behavior of $\mu_{A,2}$ in the second-half of 2020, which was a particularly dry period compared to previous years in La R\'eunion. Note also that when mechanical control is considered, it will impact $\mu_{A,2}$: for instance, thanks to Formula \eqref{muA2}, $40\%$ of Mechanical control increases $\mu_{A,2}$ by $66.7\%$.

\begin{figure}[!ht]
\centering
\begin{center}
\includegraphics[width=1.0 \linewidth]{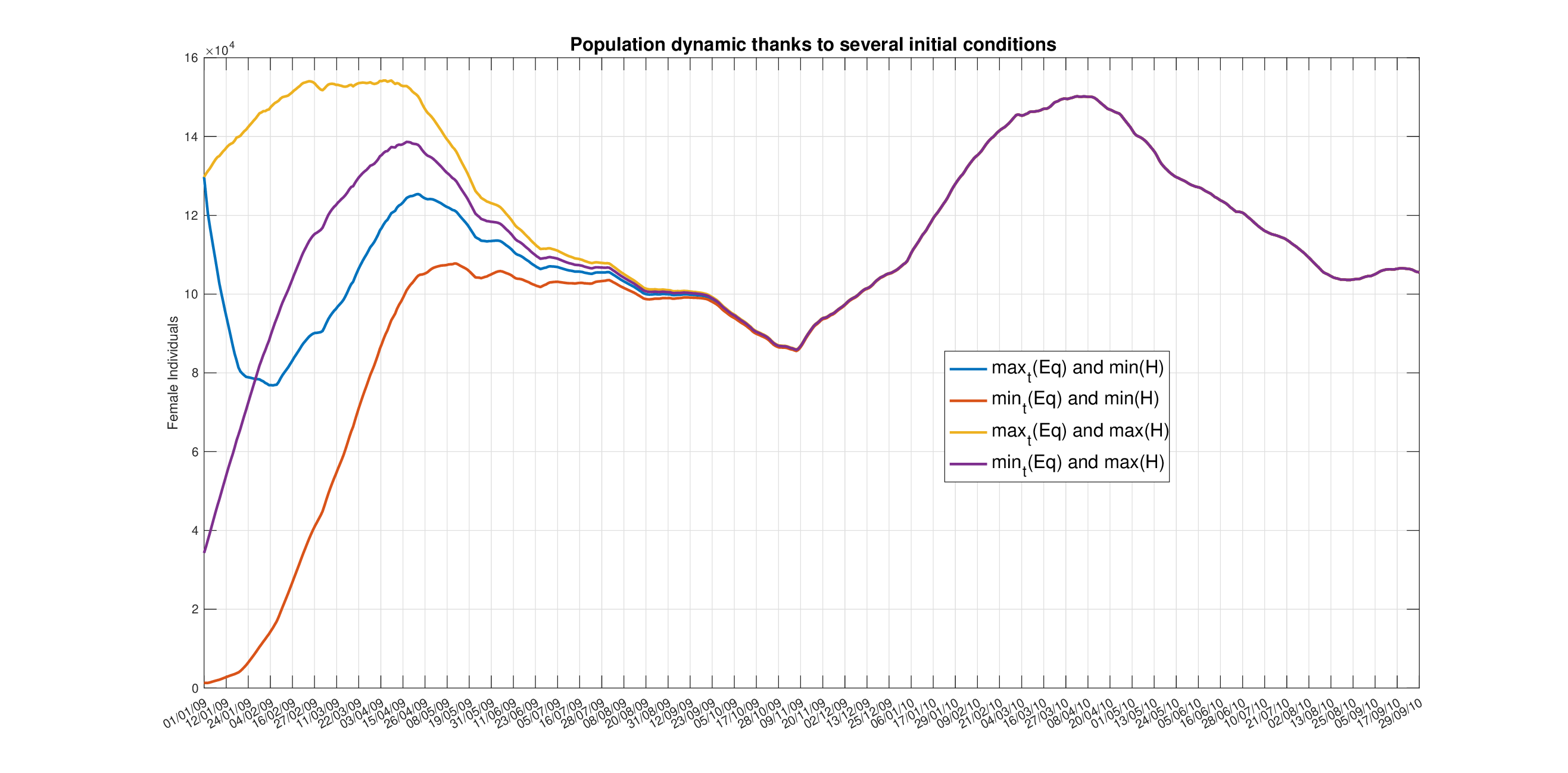} 
\end{center}
\caption{Temperature and rainfall dependent model - Simulations of the mosquito dynamics with several initial conditions, thanks to the initial rainfall data and the initial size of the mosquito population}
\label{fig:init H}
\end{figure}

The choice of the initial condition for $H$ and the initial conditions for the mosquitoes will impact respectively the initialization of the carrying capacity and the dynamic of the system for at most $12$ months (see Fig. \ref{fig:init H}, page \pageref{fig:init H}, where we have considered extreme initial values for $H(0)$ and for the  population). 
That is why, in order to have a reliable estimate of the wild mosquito population at the beginning of the control, i.e. the 1st of Sept. 2010, we start the simulations at the beginning of January 2009.

\begin{figure}[!h]
\centering
\begin{center}
\includegraphics[width=1.0 \linewidth]{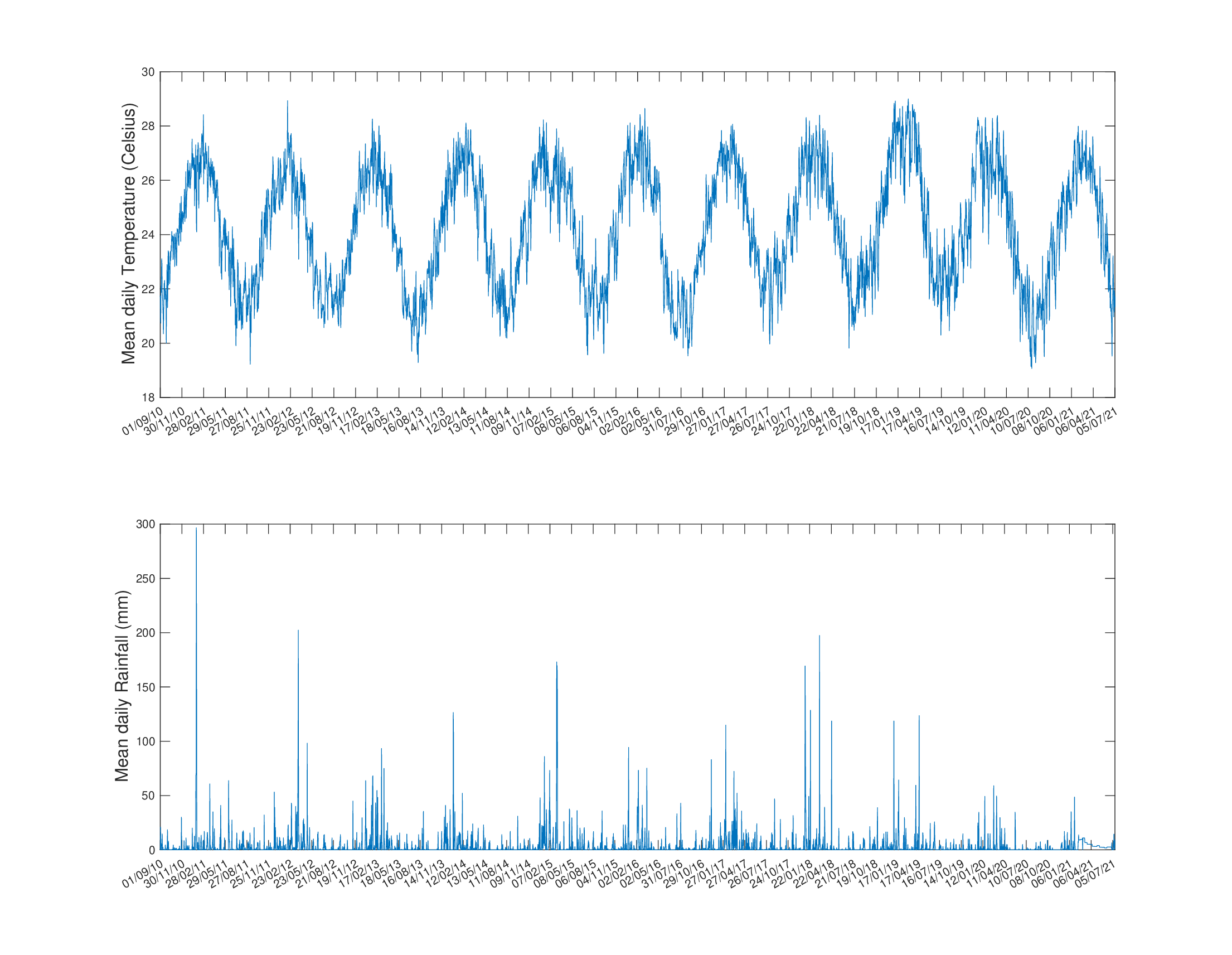}
\includegraphics[width=1.0 \linewidth]{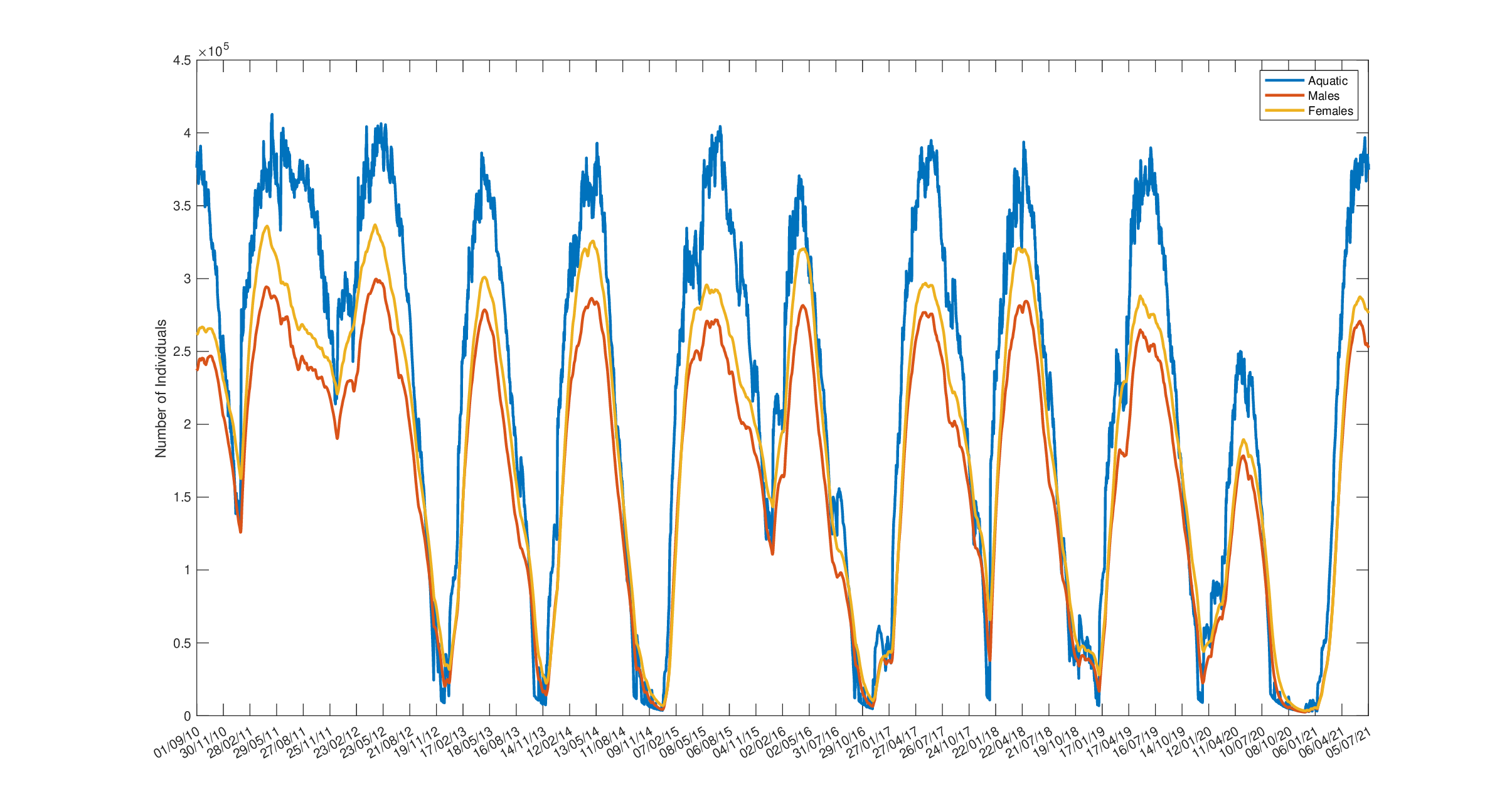} 
\end{center}
\caption{Temperature and rainfall dependent model: (a) Mean Daily Temperature (Celsius) - (b) Mean daily rainfall (mm) - (c) Mosquito population dynamics from the 1st of September 2010 until the 14th of July 2021}     \label{fig:2}
\end{figure}

Using the parameters values and the estimates of $\mu_{A,2}$ given in Fig. \ref{fig:1}, page \pageref{fig:1}, we derive the dynamic of the mosquito population in Duparc without release: see Fig. \ref{fig:2}(c), page \pageref{fig:2}. \rewieverTwo{In Fig. \ref{fig:2}(a)-(b), page \pageref{fig:2}, we also show the time evolution if the mean daily Temperature and the mean daily rainfall (mm).} As expected, periods, where the rainfall is low, imply a rapid decay of the population size, leading to an "almost" oscillatory behavior. This result is confirmed by Mark-Release-Recapture experiments derived in Duparc \cite{Legoff2019}, where the ratio between the mosquito density/ha within the dry period and the mosquito density/ha within the wet period is a factor $10$, which highlights the importance to consider a temporal dynamic in the parameters.

\subsection{Global sensitivity analysis}
\rewieverTwo{Before starting the simulations, we provide a global sensitivity analysis (GSA) using Latin hypercube sampling (LHS) and Partial rank correlation coefficient(PRCC), within the range of values for the temperature and rainfall given in Fig. \ref{fig:2}(a)-(b), page \pageref{fig:2}. The main objective of this GSA is to identify key parameters that most drive the variables and the dynamics of our system. Briefly, the objective of the LHS-PRCC sensitivity analysis is to identify key parameters whose uncertainties contribute to prediction imprecision and to rank these parameters by their importance in contributing to this imprecision. We refer the interested reader to \cite{Marino2008,Afsar2021} for further information about LHS-PRCC analysis. The forthcoming LHS-PRCC figures are made using R\cite{R} and R-studio \cite{Rstudio}.
\begin{figure}[!h]
\centering
\begin{center}
\includegraphics[width=0.7 \linewidth]{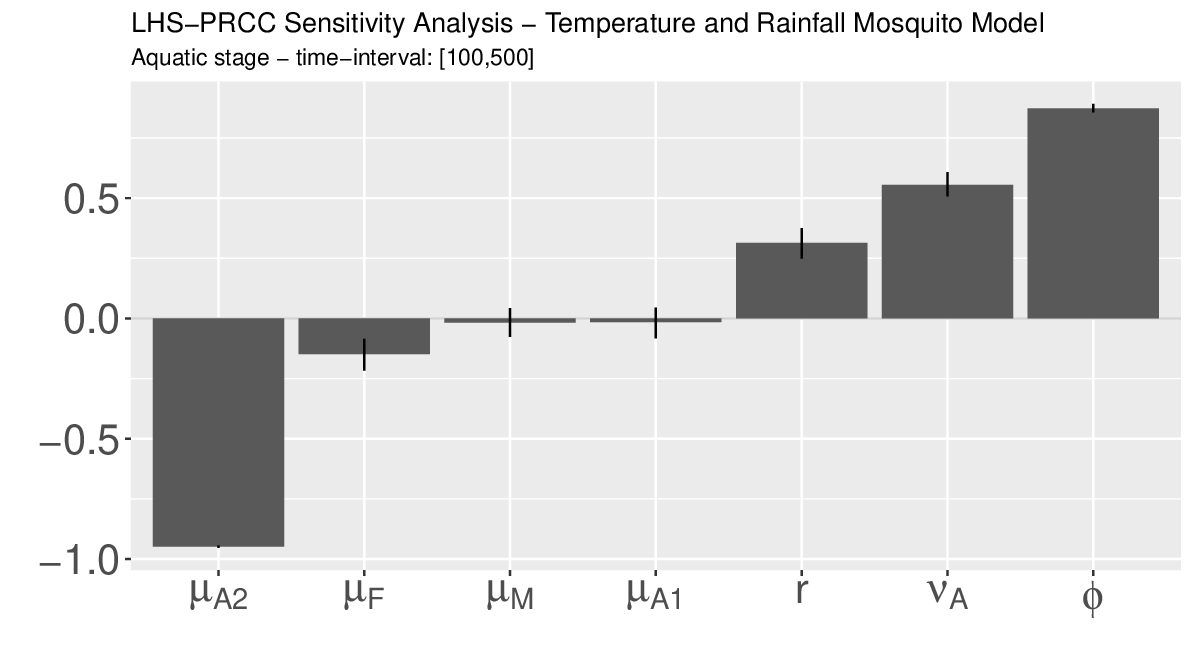}
\includegraphics[width=0.7 \linewidth]{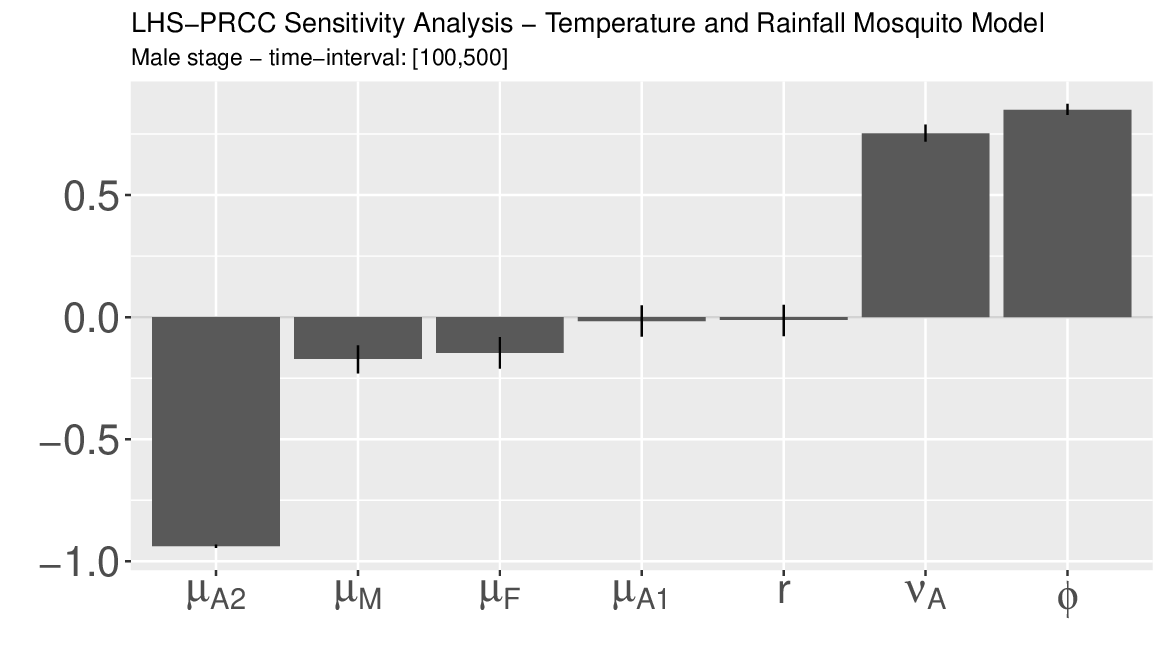} 
\includegraphics[width=0.7 \linewidth]{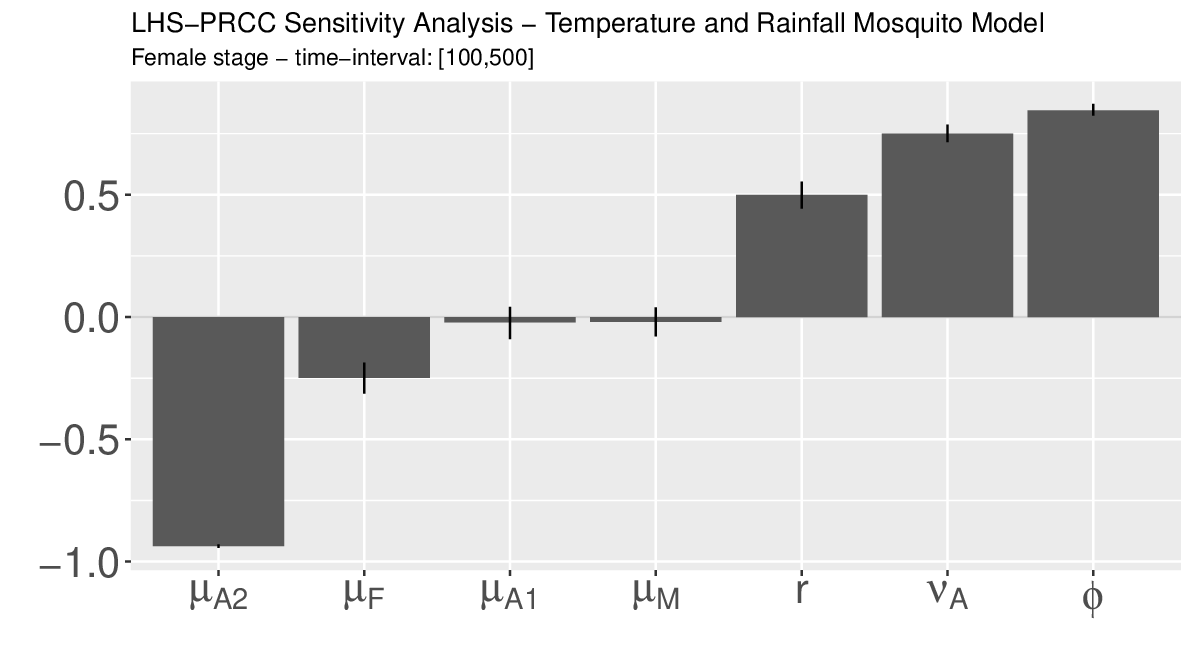} 
\end{center}
\caption{LHS-PRCC analysis of model \eqref{eq:no SIT}} \label{fig:LHS-PRCC-model1}
\end{figure}
}

\rewieverTwo{In Fig. \ref{fig:LHS-PRCC-model1}, page \pageref{fig:LHS-PRCC-model1}, the main sensitive parameters for all variables of model \eqref{eq:no SIT} are $\mu_{A,2}$, $\Phi$ and $\nu_A$, as expected. An important outcome is that since $\mu_{A,2}$ is strongly related to the temperature and rainfall,  mechanical control, i.e. the removal of breeding sites, is more than useful to increase $\mu_{A,2}$, and, thus, to impact (negatively) all compartments. This fact is also confirmed in the next LHS-PRCC analysis of the SIT model with continuous nand periodic releases: see Fig. \ref{fig:LHS-PRCC-modelTIS}, page \pageref{fig:LHS-PRCC-modelTIS}, and \ref{fig:LHS-PRCC-modelTIS_Per}, page \pageref{fig:LHS-PRCC-modelTIS_Per}. Also, as expected $\Lambda_S$ and $\beta$ are sensitive parameters in the SIT models. Contrary to the initial entomological model, the parameter $\mu_F$ has become more sensitive in the SIT models. All these informations show that a combined effort mixing SIT and mechanical control will have a strong impact on the wild population control. Last, if an additional control tool against females is able to increase $\mu_F$, using, for instance, killing traps, the overall control will be even more efficient. Last but not least, for $\tau$ varying from $1$ (daily releases) to $14$, its impact on the periodic SIT-model is negligible: what matters are the size of the releases and the competition parameters.
\begin{figure}[h]
\centering
\begin{center}
\includegraphics[width=0.7 \linewidth]{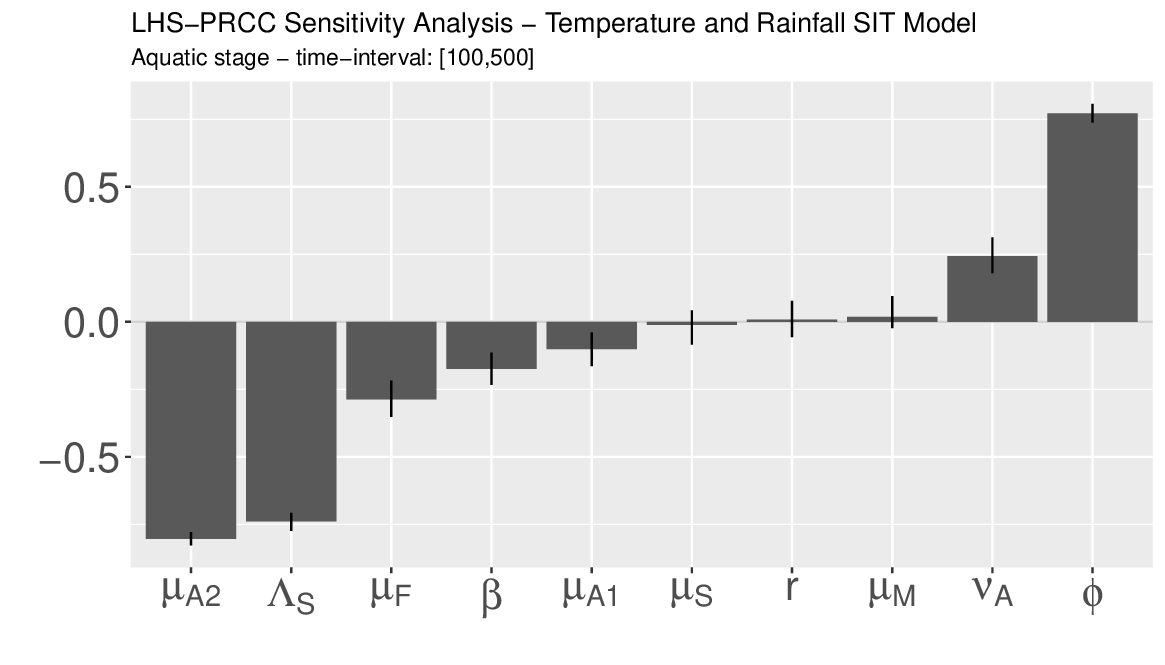}
\includegraphics[width=0.7 \linewidth]{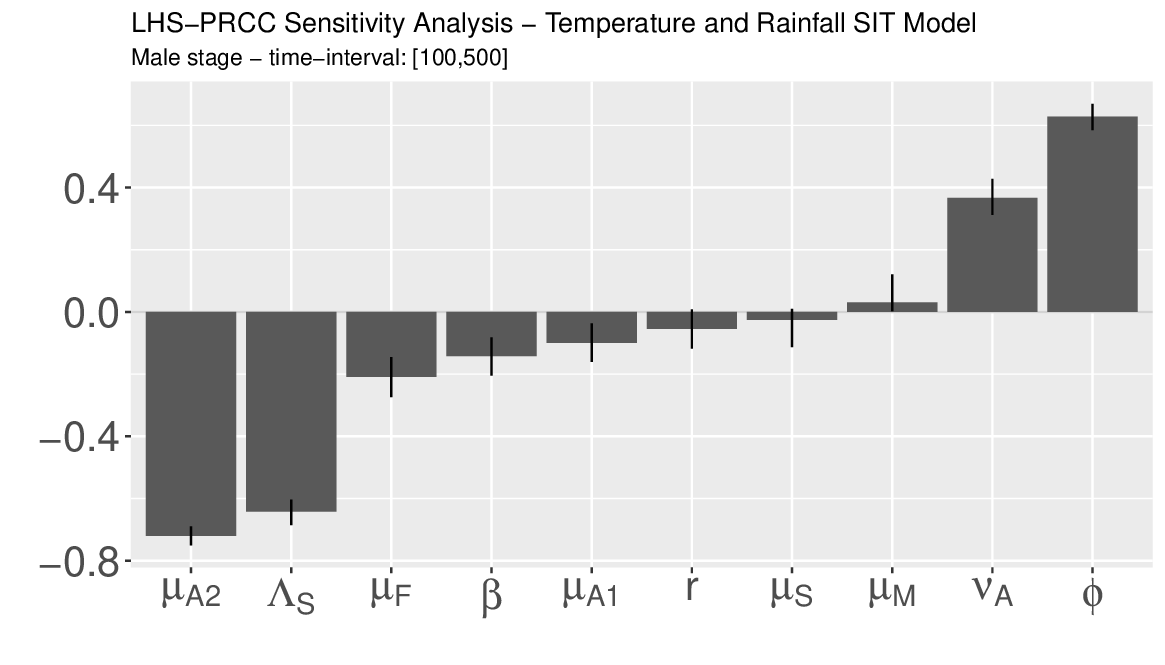} 
\includegraphics[width=0.7 \linewidth]{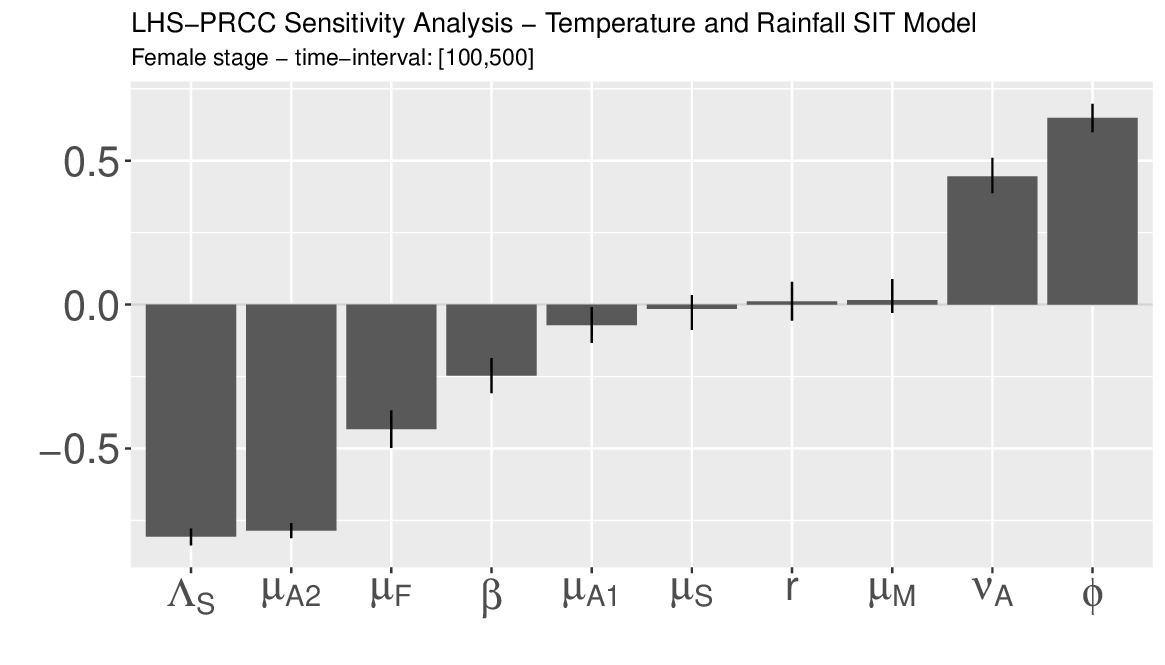} 
\end{center}
\caption{LHS-PRCC analysis of the SIT-model \eqref{TIS_mod}}
\label{fig:LHS-PRCC-modelTIS}
\end{figure}
\begin{figure}[h]
\centering
\begin{center}
\includegraphics[width=0.7 \linewidth]{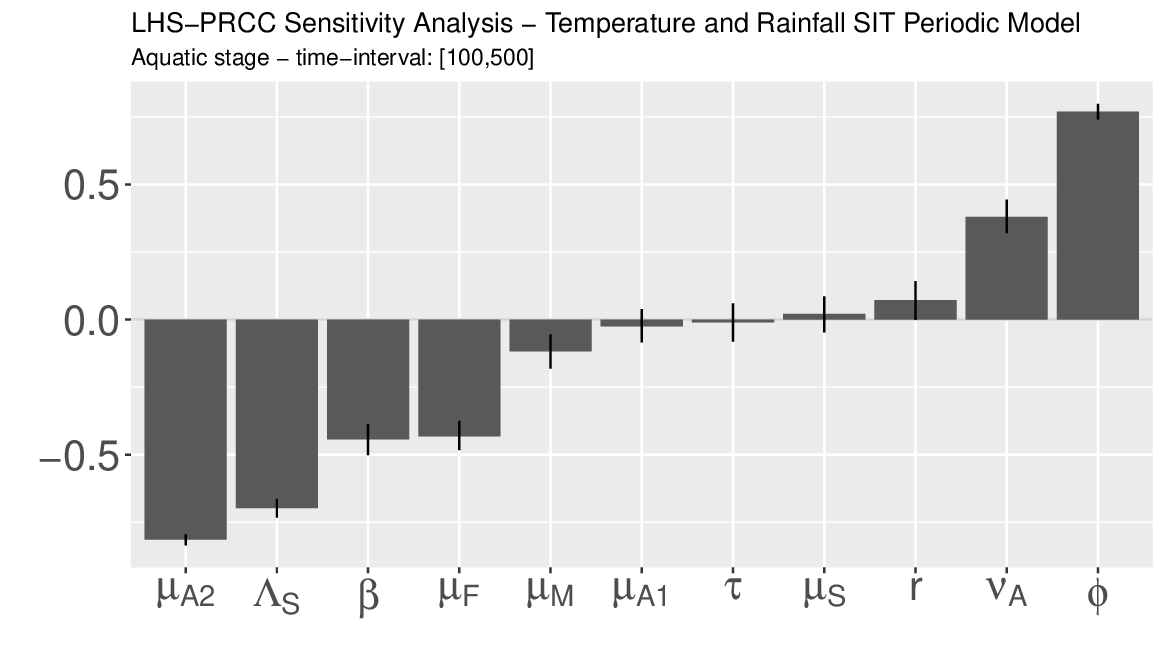}
\includegraphics[width=0.7 \linewidth]{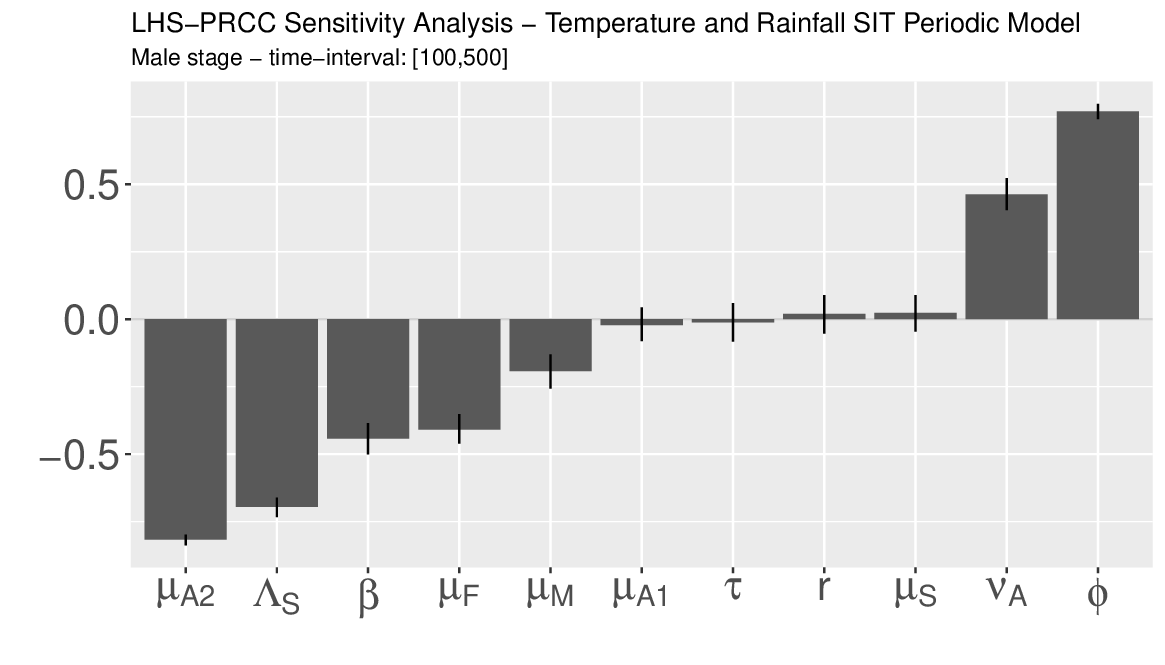} 
\includegraphics[width=0.7 \linewidth]{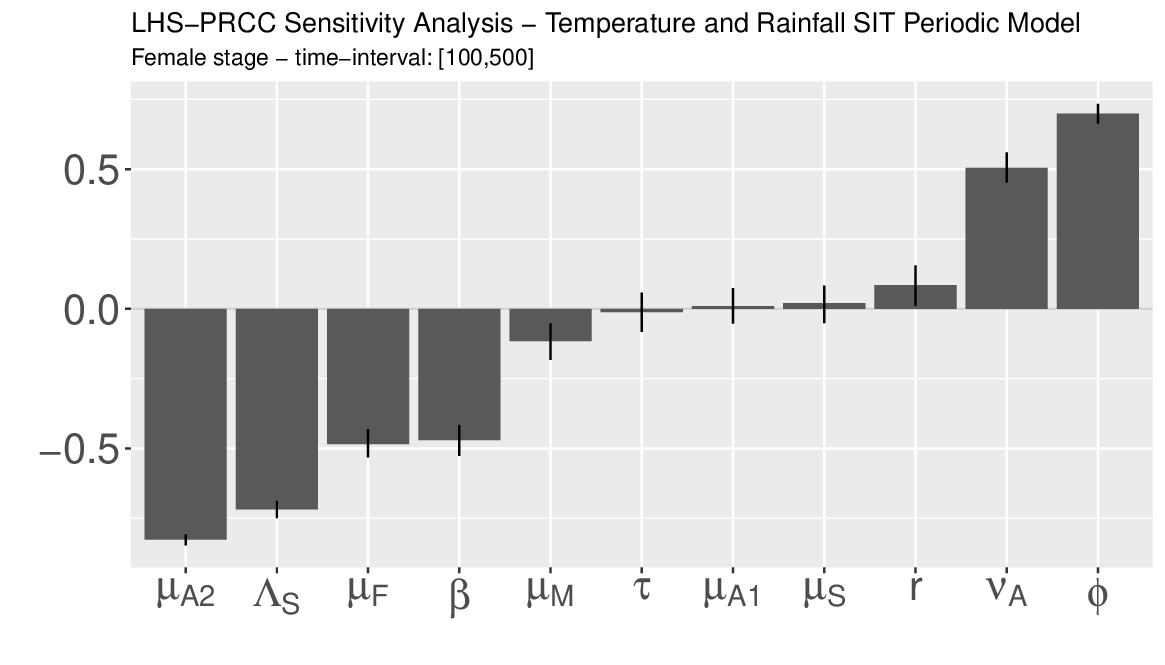} 
\end{center}
\caption{LHS-PRCC analysis of the SIT-model \eqref{TIS_mod} with periodic releases}
\label{fig:LHS-PRCC-modelTIS_Per}
\end{figure}
}
\subsection{Reducing mosquito nuisance}

As studied in \cite{AnguelovTIS2020}, we consider the massive-small strategy, which consists first of periodic massive releases, such that $\tau \Lambda_{massive}=20\times 6 000$  or $20\times 12 000$ individuals, until the wild population has become lower than $\bf{E_{1,min}}(\Lambda_{small})$, where here $\tau\Lambda_{small}=20\times 100$ individuals, such that 
$$\bf{E_{1,min}}(\Lambda_{small})\approx \left\{\begin{array}{ll}
(37.0, 29.8, 0.511)&\mbox{ if }\varepsilon=0.0,\\
(17.97, 14.47 ,0.24)&\mbox{ if }\varepsilon=0.006,\\
(1.38, 0.350, 0.0059)&\mbox{ if }\varepsilon=0.012.
\end{array}\right.$$
\yves{The previous values have been estimated thanks to the equilibria computed in Appendix B, page \pageref{AppendixB}.}
\rewieverTwo{We notice that,} when the residual fertility increases, then the size of the box $[0,\bf{E_{1,min}}(\Lambda_{small})]$ becomes smaller: SIT introduces a strong Allee effect while the residual fertility weakens it. Thus, in the forthcoming simulations, it will not be surprising to find large duration time values when the residual fertility is large. Once $[0,\bf{E_{1,min}}(\Lambda_{small})]$ is reached, the periodic releases continue at rate $\Lambda_{small}$. Using the numerical simulations, we are able to estimate the date where the system can switch from massive releases to small releases and thus evaluate the duration of the (very) massive releases.


According to the time variation of the parameters, the basic offspring number $\mathcal{N}$ will vary between $29.7$ and $85.7$ with a mean of $48.41$ over the considered period. 
Thus even if Formula \eqref{condeps} is not verified at each time, it is satisfied in mean, i.e.
\begin{equation}\label{eq:1/N}
\varepsilon<\frac{1}{T}\int_0^T\frac{1}{\mathcal{N}(t)}dt=0.0222.
\end{equation}



An important issue in SIT control is to estimate the duration of the massive releases. Thanks to contrasted environmental conditions, within the period $[2010-2020]$, we show that the minimal time to decay the wild population under a certain threshold can greatly vary.

In general, in the literature, many SIT simulations are done with constant parameters, except, for instance, in \cite{Dufourd2013,douchet2021comparing}. If we consider the mean value for each parameter over the whole period, we derive $\N \approx 49.3$ which is a quite large basic offspring number, but not surprising within a tropical context. 

\rewieverTwo{Before starting the simulations, we would like to emphasize that our meteorological data last from 01/01/2009 and 14/07/2021. Our simulations intend to derive an estimate of SIT treatment duration within this interval when the SIT control starts 01/09/2010. Thus, the reader does not have to be surprised to see that our results curves do not end at the same time: this is simply because the simulations were not able to enter $[0,\bf{E_{1,min}}(\Lambda_{small})]$ before the 14th of July 2021. Thus the simulation stops.}

\yves{In the forthcoming simulations of model \eqref{TIS_mod}-\eqref{impulsive}, we will consider four cases to illustrate the importance of considering all (or not)  environmental parameters:
}
\begin{enumerate}
\item \yves{Model 1:} model \eqref{TIS_mod}-\eqref{impulsive} with temperature and rainfall dependent parameters.
\item \yves{Model 2:} model \eqref{TIS_mod}-\eqref{impulsive} with temperature-dependent parameters only.
\item \yves{Model 3:} model \eqref{TIS_mod}-\eqref{impulsive} with a constant average temperature and rainfall-dependent parameters.
\item \yves{Model 4:} model \eqref{TIS_mod}-\eqref{impulsive} with average parameters values estimated from September 2010 to mid-July 2021.
\end{enumerate}

\begin{figure}[!ht]
\centering
\includegraphics[width=1.0 \linewidth]{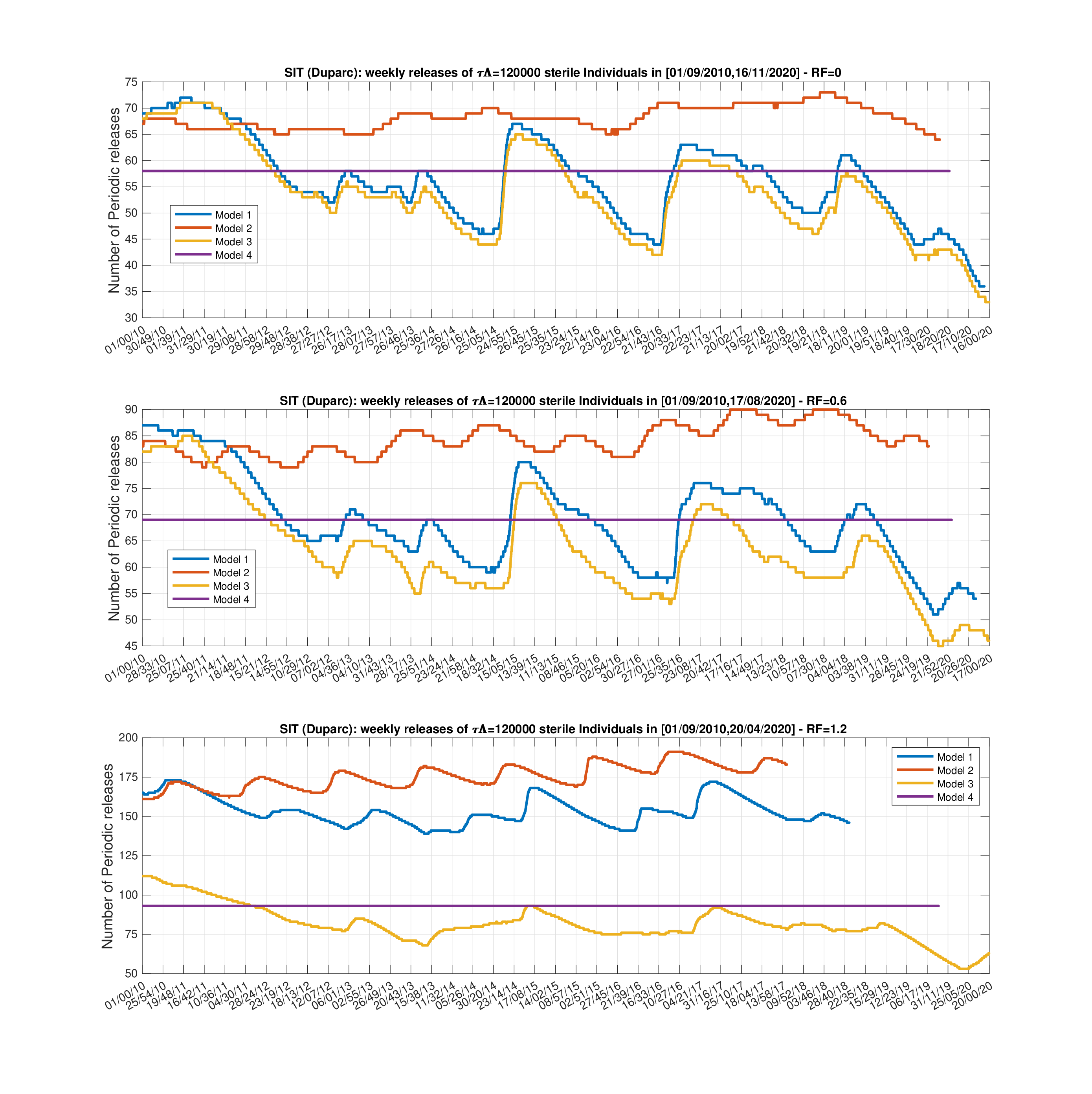} 
\caption{Comparison of the amount of weekly SIT releases calculated using the four models: (a) $0\%$ of residual fertility; (b) $0.6\%$ of residual fertility; (c) $1.2\%$ of residual fertility without mechanical control - The weekly release rate is $\tau\Lambda= 120000$ Individuals}
\label{comparison_nuisance_models}
\end{figure}

\begin{figure}[!ht]
\centering
\includegraphics[width=1.0 \linewidth]{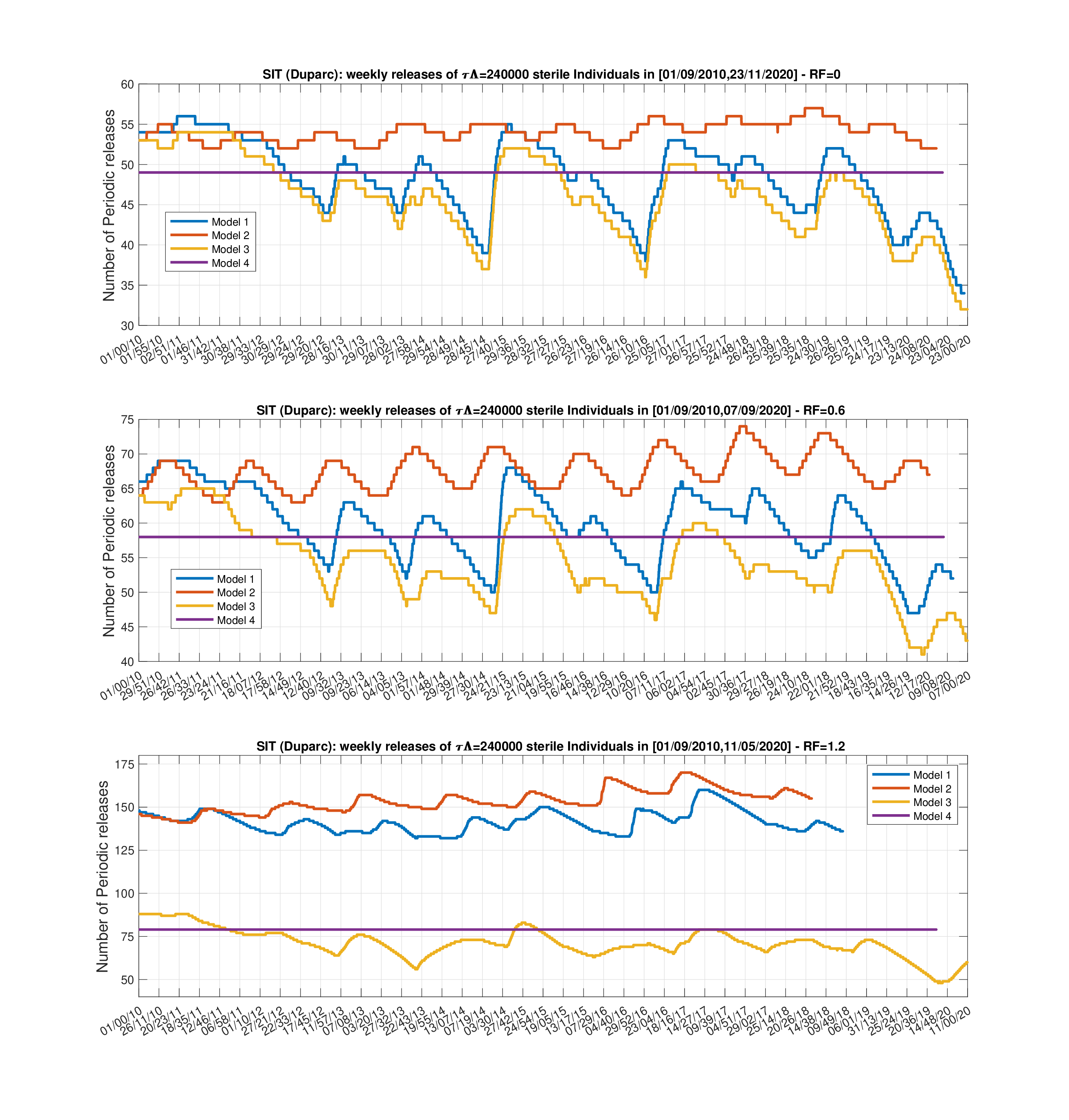} 
\caption{Comparison of the amount of weekly SIT releases calculated using the four models: (a) $0\%$ of residual fertility; (b) $0.6\%$ of residual fertility; (c) $1.2\%$ of residual fertility without mechanical control - The weekly release rate is $\tau\Lambda= 240000$ Individuals.}
\label{comparison_nuisance_models240000}
\end{figure}

\yves{We first provide Fig. \ref{comparison_nuisance_models}, page \pageref{comparison_nuisance_models}, where simulations of the four models are given for different values of RF, $0\%$, $0.6\%$ and $1.2\%$, with $\tau \Lambda=120000$ individuals, and without mechanical control. It is easy to check that the greater RF, the larger the amount of releases. However, the increase is reasonable when $RF=0.6\%$, it is really important when $RF=1.2\%$. It is interesting to notice that model 1 and model 3 behave similarly more or less: this is due to the fact that rainfall is taken into account in parameter $\mu_{A,2}$. Thus, for RF equal to $0\%$ and $0.6\%$, model 3 is close to the variations of model 1. In contrary, when $RF=1.2\%$, then it seems that the temperature-dependent parameters are taken over since model 1 and model 2 are now close, while model 3 provides very different and lower estimates. Model 4 provides a constant value for the amount of releases: while this average value seems to be realistic when residual fertility is small, this is no more the case when residual fertility is large: see Figs. \ref{comparison_nuisance_models} and \ref{comparison_nuisance_models240000}. Model 2, where the parameters are only temperature dependent, provide the worst estimates, except in the years where rainfall was abundant, i.e. 2010 and 2011. However, along the decade, except when RF is equal to $1.2\%$, model 2 is the worst one. In addition, the main disappointment with model 2, with temperature-dependent parameters only, comes from the fact that the inter-annual periodic behavior, for which the best period to start SIT would (roughly) be between May and September. This is exactly the same conclusion reached by entomologists.}

\yves{In Fig. \ref{comparison_nuisance_models240000}, page \pageref{comparison_nuisance_models240000}, we consider the same simulations than in Fig. \ref{comparison_nuisance_models}, page \pageref{comparison_nuisance_models}, but with the release of $\tau \Lambda=240 000$ individuals per week. The results, in terms of behavior, are almost the same with some improvements in terms of the releases amount. However, thanks to the sterile males production effort, the gain are very little. This shows that increasing the release rate does not necessarily improve the duration of massive releases.}

\begin{figure}[!ht]
\centering
\includegraphics[width=1.0 \linewidth]{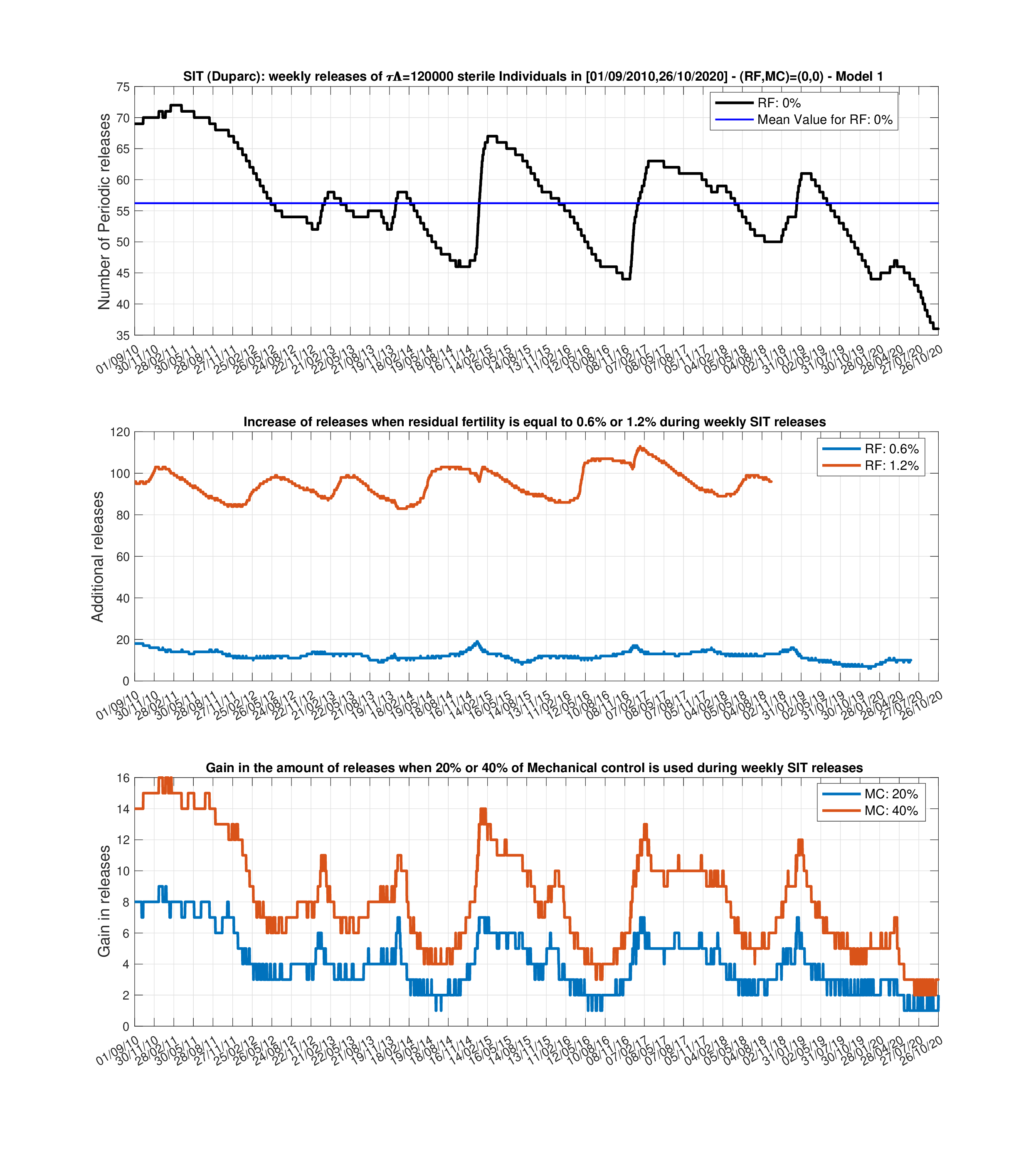} 
\caption{Simulations of Model 1 with $\tau\Lambda= 120000$ Individuals: (a) Variation of the amount of SIT releases; (b); Impact in terms of additional releases when $0.6\%$ and $1.2\%$ residual fertility occurs; (b) Gain in the releases when $20\%$ and $40\%$ of Mechanical control is used}
\label{model1_variation_120000}
\end{figure}

\yves{The death rate of mosquitoes is mainly linked to the temperature, such that we can intuitively think that the best period to act is during the Austral winter, i.e. from the end of June to the end of September. However, According to Figs. \ref{comparison_nuisance_models}, \ref{comparison_nuisance_models240000}, and \ref{model1_variation_120000}, page \pageref{comparison_nuisance_models}, \pageref{comparison_nuisance_models240000} and \pageref{model1_variation_120000}, and following model 1 (and model 3), the best intervention period, i.e. the SIT starting-time, lasts, in general, from July to December.} \rewieverTwo{Indeed, as seen in these figures, in particular in Fig. \ref{model1_variation_120000}, page \pageref{model1_variation_120000}, the amount of releases oscillates with various amplitudes but with its minimal value always reached around the last two weeks of November, except in 2020 that was a particularly dry year in La R\'eunion. After November, we usually enter in the rainy season, and that is why the amount of releases increases fast. In fact, each year, from June-July to December, we observe that the amount of releases needed to reach elimination reaches its lowest value and is, more or less, decreasing from June to December.} Thus, the window to start SIT control is larger than expected and is not necessarily reduced to the Austral winter period  (from June to September) but can also include the Spring period (from October to December).


\yves{We now focus on simulations of model 1, with temperature and rainfall dependent parameters, taking into account the impact of residual fertility and Mechanical control: see Fig. \ref{model1_variation_120000}, page \pageref{model1_variation_120000}. In Fig. \ref{model1_variation_120000}(a), the straight line represents the average value of the numerical simulations. Thanks to Fig. \ref{comparison_nuisance_models}, we can see that that this value is 4 weeks larger than the value obtained by model 4, where all parameters are supposed to be constant. However, around this average value, we can see that the results are very contrasted, mainly related to the average rainfall. In Fig. In Fig. \ref{model1_variation_120000}(b), we show the impact of residual fertility: while $0.6\%$ of residual fertility "only" increase the amount releases up to 20 weeks, a residual fertility of $1.2\%$ increases drastically the amount of the releases, between $80$ and $110$ weeks. Finally, in Fig. \ref{model1_variation_120000}(c), we show the benefit of Mechanical control, i.e. reducing the breeding sites. However, we recover the fact that mechanical control is more beneficial within rainy periods, when the mosquito population is large, than in dry periods: compare the years 2010-2011 with the years 2016, 2019, and 2020, in Fig. \ref{model1_variation_120000}, page \pageref{model1_variation_120000}. This is quite obvious but this is important to highlight this fact in terms of control strategies and cost.
}

\yves{From the practical point of view, since we are not yet capable of predicting weather data accurately, it is important to have an initial estimate for the duration of the massive treatment and then refine this estimate while the process is in progress. Thus the "mean constant" simulations, obtained with model 4, provide reasonable values, at least for small residual fertility, that can help to evaluate an average duration of the massive releases with and without mechanical control.} However, this approach does not provide any information on the best period to start the SIT treatment to reduce the duration of the massive releases and, thus, optimize the sterile males production.

Based on the average value given by model 4, we can estimate, roughly at least, the duration and the amount production of sterile males needed for the massive releases to reduce the nuisance. In general, the amount of sterile males to produce is considerable. However, all along the treatment, the size of the massive releases can be adapted to the wild male population, using Mark-Release-Recapture experiments, as well as the periodicity of the releases, like in \cite{Bliman2019,Aronna2020}. 

Last, if the residual fertility is too large, say to $2\%$, then it is impossible to reduce (in a reasonable amount of time) the wild population under the given threshold, such that the wild population can be controlled with small releases of $\tau \Lambda_M=20 \times 100$ individuals. The recommendation is thus to improve the sterilization process to have a residual fertility as small as possible, at least, less than $0.6\%$, for instance.

Our simulations took place in the context of nuisance reduction, i.e. to reduce the wild mosquito population in order to reduce the number of bites. It is not always necessary to reach this objective. In particular, in a tropical context, where people are used to mosquitoes. There, the most important goal is to reduce the epidemiological risk.
%
%
\subsection{Reducing the epidemiological risk}

In the previous simulations, we derive numerics to lower the mosquito population under a given release threshold for sterile males, for instance, $100$ Ind/ha/week. As explained, in a tropical context, another option is to reduce the population in order to reduce the epidemiological risk, to prevent the risk of an epidemic. Since dengue is often circulating in La R\'eunion, we can couple our entomological model with a dengue model, like the one developed in \cite{DUMONT2022}. According to the epidemiological model developed in \cite{DUMONT2022}, and recalled in Appendix C, page \pageref{AppendixC}, we derive the following formula for the SIT basic reproduction number
\begin{equation}\label{R0-SIT}
 \mathcal{R}^2_{0,SIT}=
\dfrac{\nu_m}{\nu_m+\mu_F}\dfrac{B\beta_{mh}}{\mu_F}\dfrac{B\beta_{hm}}{\eta_h+\mu_h}\dfrac{F_S^*}{N_h}, 
\end{equation}
where $F_S^*$ is the amount of susceptible adult females at the Disease Free Equilibrium (DFE) equilibrium, estimated according to the value taken by $\tau\Lambda_{Massive}$. In \cite{DUMONT2022}, the authors showed that $F_S$ is either strictly positive or equal to zero, depending on the amount of sterile males released. The positive parameters $\mu_h$ and $1/\nu_h$ represent respectively the average human  mortality rate and the average viremic period. Since no disease-induced mortality is considered, the total human population is supposed to be constant and equal to $N_h$. The average rate of mosquito bites per individual is denoted $B>0$, and $\beta_{mh}>0$ ($\beta_{hm}>0$) is the probability of dengue transmission from an infected female mosquito (human) to a susceptible human (mosquito) during such an event. The positive parameter $\nu_m>0$ is the extrinsic incubation rate (EIR). However, it is well known that the previous \rewieverTwo{epidemiological parameters, $B$,} $\nu_m$, $\beta_{hm}$, and $\beta_{mh}$, are Temperature-dependent. Using results from \cite{Xiao2014}, we consider the following transmission
probability for an infected \textit{Aedes albopictus} to transmit DENV-2, using a Lactin-1 function
\[
\beta_{mh}(T)=
\exp(\alpha\times T)-\exp\left(\alpha \times T_{max}-\dfrac{(T_{max}-T)}{\delta_{T}})\right),
\]
with $\alpha=0.20404$, $T_{max}=37.354$, and $\delta_{T}=4.89694$.
The probability of transmission from humans to mosquitoes is negligible for "low" temperatures, increases linearly to one at a maximum temperature and remains at one for higher temperatures. Thus, following \cite{Naven2019}, we consider 
\[
\beta_{hm}(T)=\dfrac{T^{7}}{T^{7}+\beta_{h}^{7}},
\]
with $\beta_{h}=18.9871$. Last but not least, the EIP (Extrinsic incubation Period) decreases according to the temperature. Using again \cite{Xiao2014}, we derive the following interpolation:
\[
\nu_m(T)=aT^{2}+bT+c,
\]
with $a=-0.001$, $b=0.0670$, and $c=-0.866$. All other (epidemiological) parameters are supposed to be constant.

\rewieverTwo{Finally, following \cite{Mordecai2017}[Table A in the supplementary material], we consider the following Bri\`ere function for the biting rate $B$:
$$
B(T)=c\times T \times \left(T-T_{\min}\right)\sqrt{T_{\max}-T},
$$
where $c=1.93e-04$, $T_{\min}=10.25$, and $T_{\max}=38.32$ are the parameters related to \textit{Aedes albopictus} mosquito. 
}

The threshold $\mathcal{R}_{0,SIT_c}^2$ is related to the long-time behavior of the system. For practical purpose, we will consider $\mathcal{R}_{eff}$, the effective reproduction number, that is defined as follow 
\begin{equation}
    \mathcal{R}_{eff}(t)=\dfrac{\nu_m(t)}{\nu_m(t)+\mu_F(t)}\dfrac{B(t)^2\beta_{mh}(t)\beta_{hm}(t)}{\mu_F(t)\left(\eta_h+\mu_h\right)}\dfrac{F_S(t)}{N_h}.
\end{equation}

In fact, assuming that $t_{DENV}$ is the time where a DENV virus starts circulating, we will estimate $\mR_{eff}$ at time $t_{DENV}$. Clearly, if $\mR_{eff}(t_{DENV})<1$ and $\mR_{0,SIT_c}^2<1$, then no epidemics will occur. In contrary, even if $\mR_{0,SIT_c}^2<1$ but $\mR_{eff}(t_{DENV})>1$ then an outbreak may occur.

In the forthcoming simulations, we will estimate the time needed to lower $\mR_{eff}$ below $0.5$ for different sizes of massive releases, i.e. $6 000$ or $12 000$ sterile males per ha, and for different residual fertility, with and without mechanical control. Thus, it suffices to find $t^*$ such that the wild female population verifies
\yves{
\begin{equation}
    F(t^*)<\dfrac{\nu_m(t^*)+\mu_F(t^*)}{\nu_m(t^*)}\dfrac{\mu_F(t^*)\left(\eta_h+\mu_h\right)}{B(t^*)^2\beta_{mh}(t^*)\beta_{hm}(t^*)}\dfrac{N_h}{2},
    \label{F_condition_reduce_risk}
\end{equation}
}
for a given $N_h$. For the numerical simulations, we will consider the parameters values given in Table 
 \ref{Table-valeurs-parametre-epidemio}, page \pageref{Table-valeurs-parametre-epidemio}.

\begin{figure}[!ht]
\centering
\includegraphics[width=1.0 \linewidth]{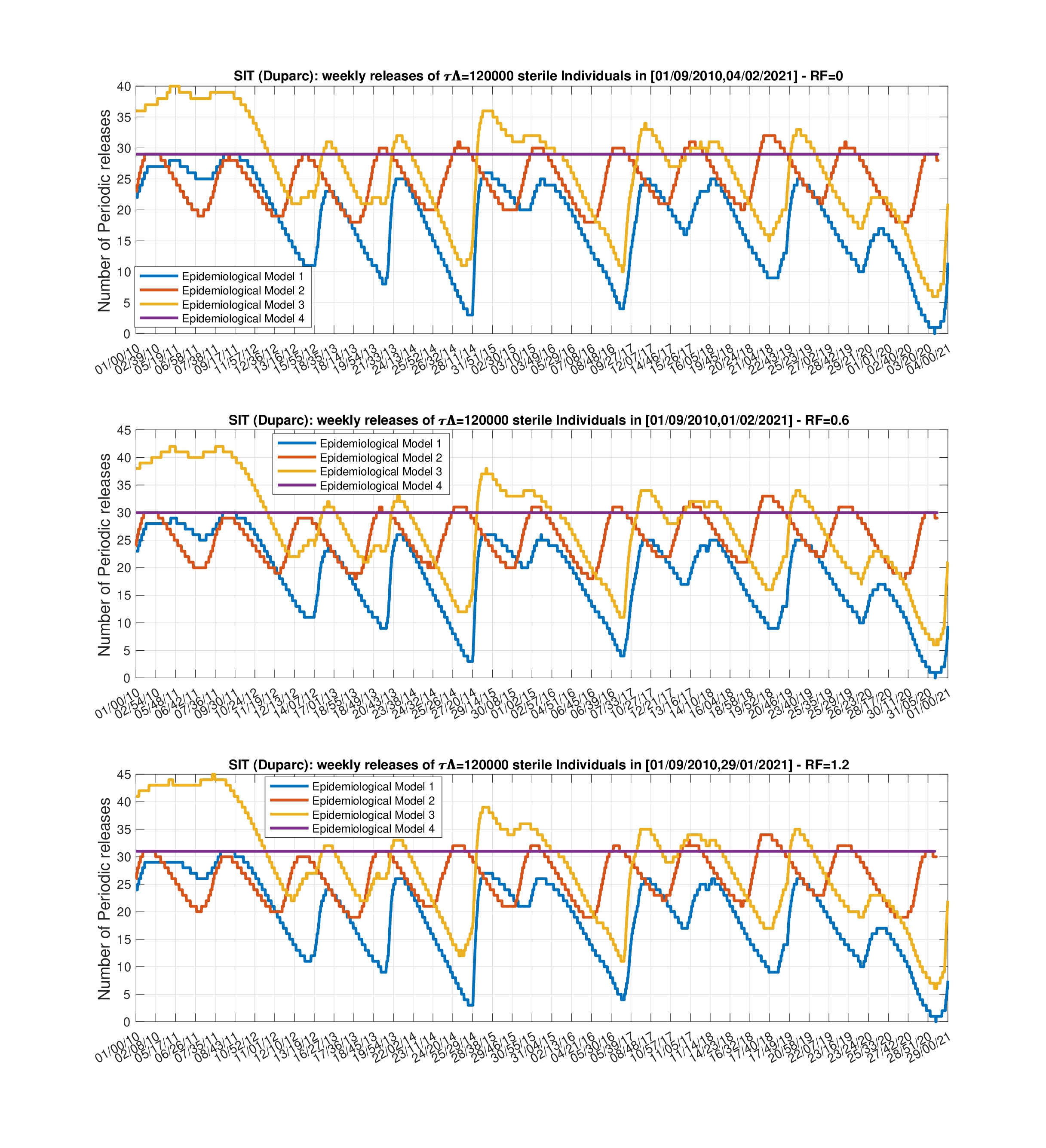} 
\caption{Reducing the epidemiological risk. Comparison of the amount of weekly SIT releases calculated using the four models: (a) $0\%$ of residual fertility; (b) $0.6\%$ of residual fertility; (c) $1.2\%$ of residual fertility without mechanical control - The weekly release rate is $\tau\Lambda= 120000$ Individuals}
\label{comparison_epidemiological_models}
\end{figure}

\begin{figure}[!ht]
\centering
\includegraphics[width=1.0 \linewidth]{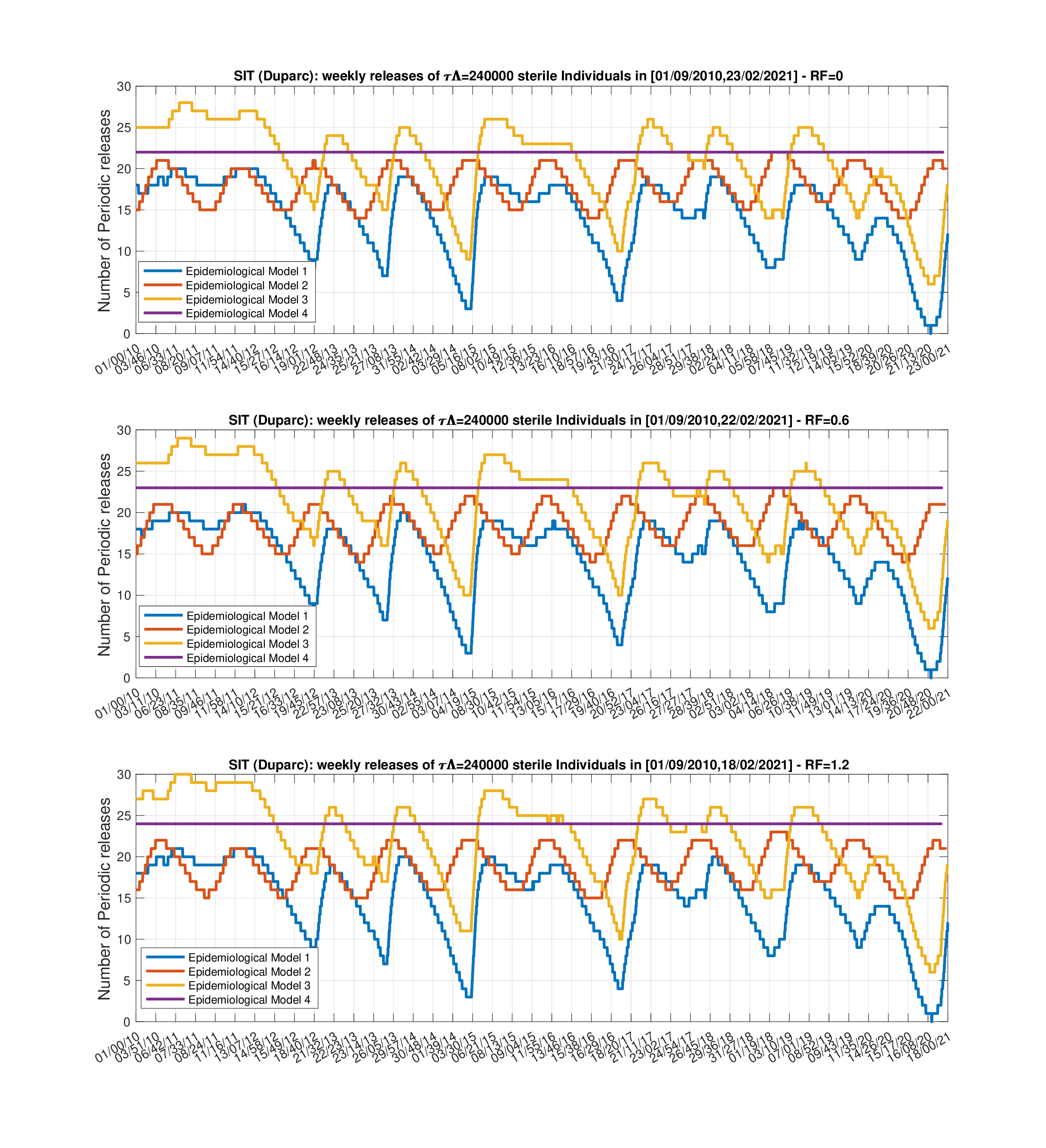} 
\caption{Reducing the epidemiological risk. Comparison of the amount of weekly SIT releases calculated using the four models: (a) $0\%$ of residual fertility; (b) $0.6\%$ of residual fertility; (c) $1.2\%$ of residual fertility without mechanical control - The weekly release rate is $\tau\Lambda= 240000$ Individuals.}
\label{comparison_epidemiological_models240000}
\end{figure}

\begin{figure}[!ht]
\centering
\includegraphics[width=1.0 \linewidth]{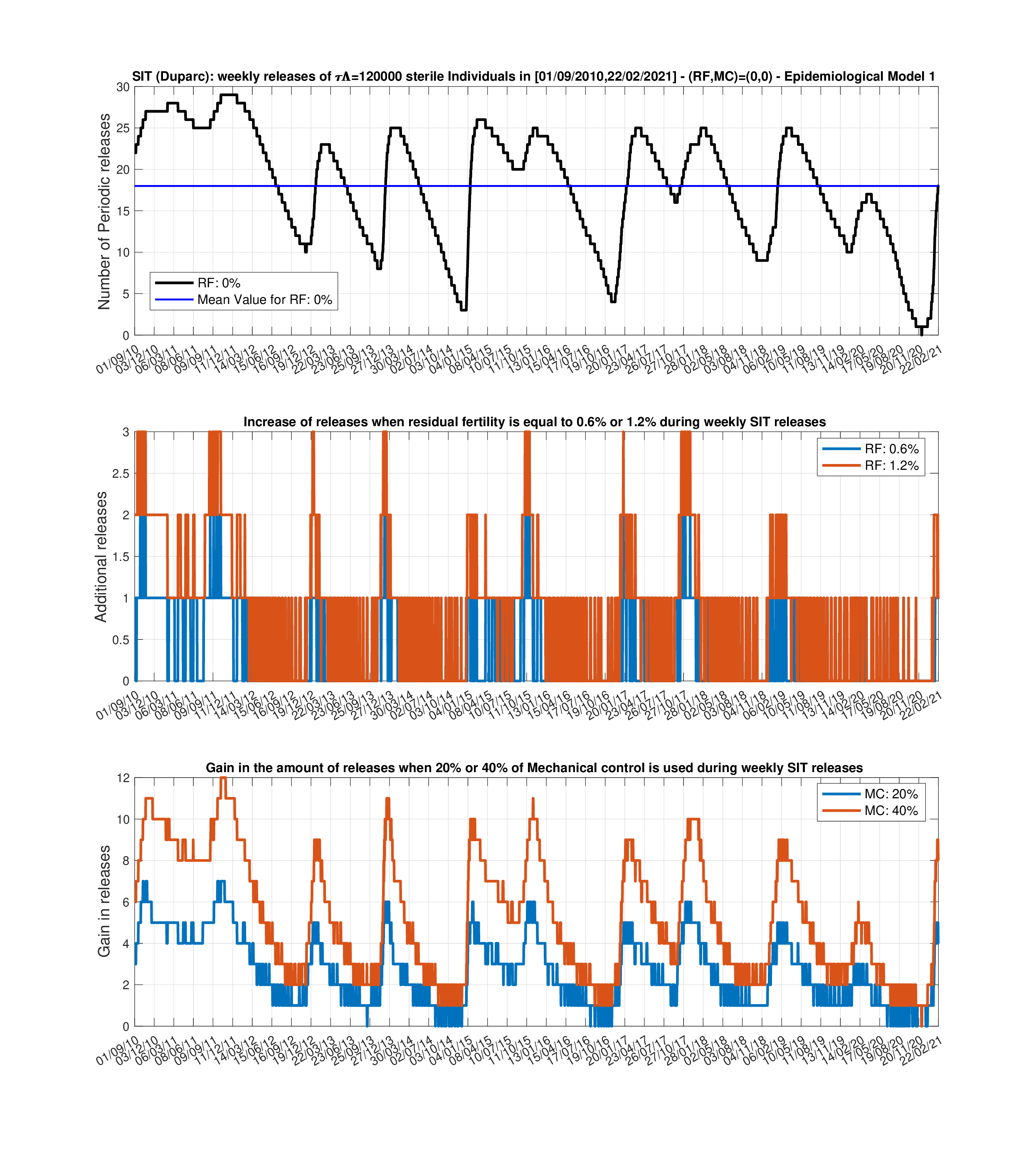} 
\caption{Reducing the epidemiological risk. Simulations of Model 1 with $\tau\Lambda= 120000$ Individuals: (a) Variation of the amount of SIT releases; (b); Impact in terms of additional releases when $0.6\%$ and $1.2\%$ residual fertility occurs; (b) Gain in the releases when $20\%$ and $40\%$ of Mechanical control is used}
\label{epidemiological_model1_variation_120000}
\end{figure}

\begin{table}[tbp]
  \begin{center}
  \caption{\bf \textit{DENV} epidemiological constant parameter values \cite{DUMONT2022}}
  \begin{tabular}{|c|c|c|c|c|}
  \hline
  Symbol &  $\mu_h$ & $\eta_h$ & $ N_h$\\
  \hline
  Value & $\dfrac{1}{365\times 78}$ & 1/7 & $2 000$ \\
  \hline
 Unit  & day$^{-1}$ & day$^{-1}$ & Ind \\
      \hline
\end{tabular}
\label{Table-valeurs-parametre-epidemio}
\end{center}
\end{table}
\yves{Assume that the weekly release rate is $\tau \Lambda= 120 000$ individuals. In Fig. \ref{comparison_epidemiological_models}, page \pageref{comparison_epidemiological_models}, for a given SIT starting date, we compute the number of releases necessary to reach $R_{eff}(t^*)<0.5$ for the four models, like in the nuisance reduction case. However, contrary to the nuisance reduction case, the results between model 1 and model 3 are quit different, while the results obtained with model 2 are more homogeneous and looks like periodic, varying between 19 and 31 weeks. With model 4, with constant average parameters, between 29 and 31 weeks are necessary to decrease the epidemiological risk, depending on residual fertility. Contrary to the nuisance reduction case, the impact of residual fertility is low: see Fig. \ref{comparison_epidemiological_models}(b)-(c). At some point this is a very good news. A special case occur: end of November 2020, where the wild initial population is so small that, for model 1, no release is necessary. If we assume $\tau \Lambda= 240 000$ individuals, we derive Fig. \ref{comparison_epidemiological_models240000}, page \pageref{comparison_epidemiological_models240000}. There is no substantial changes thanks to the case $\tau \Lambda= 120 000$, only a gain in the amount of releases, for instance $7$ weeks for model 4, and only in periods where the wild population is very large (2010, 2011, for instance). Again, the impact of the residual fertility is negligible. Last, we focus on simulations derived with model 1, when residual fertility and mechanical control occurs: see Fig. \ref{epidemiological_model1_variation_120000}, page \pageref{epidemiological_model1_variation_120000}. }

\begin{figure}[!ht]
\centering
\includegraphics[width=1.0 \linewidth]{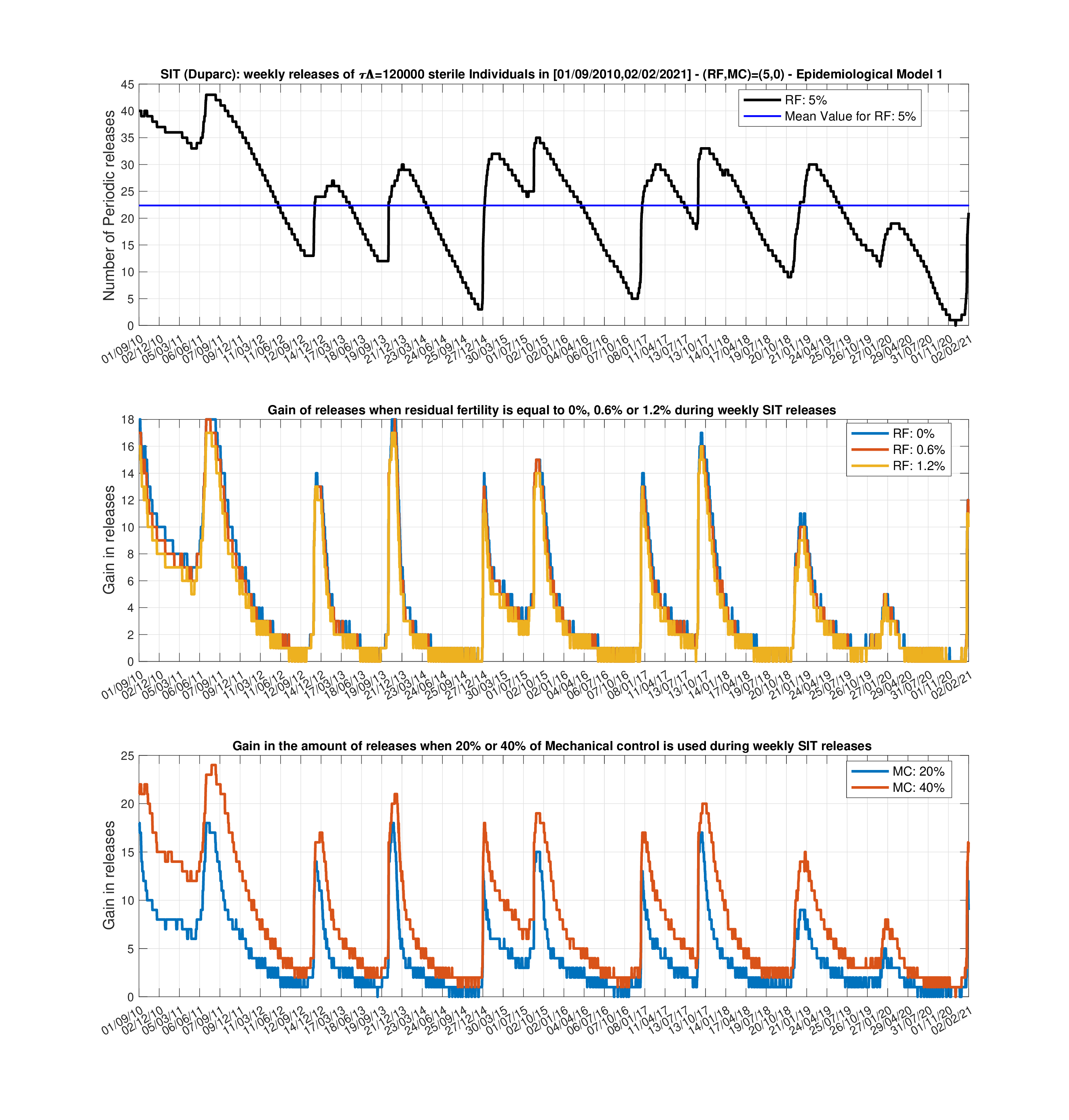} 
\caption{Reducing the epidemiological risk. Simulations of Model 1 with a residual fertility of $5\%$ and $\tau\Lambda= 120000$ Individuals: (a) Variation of the amount of SIT releases; (b); Gain in the releases when residual fertility is $0\%$, $0.6\%$ or $1.2\%$; (b) Gain in the releases when $20\%$ and $40\%$ of Mechanical control is used}
\label{epidemiological_model1_variation_120000_RF5}
\end{figure}

\yves{As expected by the previous simulations, there is no impact of the residual fertility on the amount of releases. In fact, we can even consider larger residual fertility, for instance $5\%$: see Fig. \ref{epidemiological_model1_variation_120000_RF5}, page \pageref{epidemiological_model1_variation_120000_RF5}. Like for the other residual fertility cases, the duration of SIT releases increases only within periods (of the year) where the wild population is large, while there is almost no impact of the residual fertility in periods where the population is small. In fact, the female population size to reach to get $\mathcal{R}_{eff}<0.5$, given in \eqref{F_condition_reduce_risk}, is so small that it is always above the lower bound given in \eqref{lowerbound}, page \pageref{lowerbound}, i.e. $F(t^*)>F_l^*$. This is very important to have in mind: if from time to time the residual fertility may change thanks to issues during the ionization process, our results show that this may no impact the objective of reducing the epidemiological risk.}

\yves{However, like in the nuisance reduction case, mechanical control may have an impact, but only in periods where the mosquito population is large, i.e. when rainfall are abundant: the amount of releases can be divided by a factor $2$. Otherwise, when the population is small (from September to December in general), SIT alone, with massive releases, could work. Thanks to Fig. \ref{epidemiological_model1_variation_120000} and Fig. \ref{epidemiological_model1_variation_120000_RF5}, the larger the residual fertility, the more useful the mechanical control during the rainy season.}

Of course, once $\mathcal{R}(t^*)<0.5$ is reached, it will be necessary to continue to release enough sterile males to maintain $\mathcal{R}(t^*)$ below $0.5$, as long as needed: \yves{the smaller the wild population after the massive releases, the smaller the releases to keep $\mathcal{R}(t^*)<0.5$.} This can be evaluated through numerical simulations.
\begin{remark}
In fact for this particular objective of reducing the epidemiological risk, and because in La R\'eunion we have some seasonality, the meaningful strategy, as already proposed in \cite[Fig.4]{Dufourd2013}, would be to consider massive releases only over the wet period to maintain the mosquito population at a level corresponding to the population size during the dry period, the Austral winter, where, in general, the epidemiological risk is low because the vector population is low, except, of course, when rainfalls occur during this period, like in 2010 and 2011. Thus, we would have a "massive and stop" SIT strategy, contrary to the nuisance reduction where a "massive and small" releases strategy seems more appropriate. 
Thus, according to the temperature and rainfall dependent model, massive releases could only occur from October-November to May-June, that is over $8$ months, with, eventually, a reduction of the size of the releases, using, for instance, a closed-loop control strategy \cite{Aronna2020,Bliman2019}.
\end{remark}


\section{Conclusion}
We \yves{have built} a minimalistic rainfall-temperature entomological model to derive the dynamics of \textit{Aedes albopictus} population in the place of Duparc (Sainte-Marie), La R\'eunion island. Since we are in a sub-tropical environment, our simulation shows that rainfall definitively plays a great role in the dynamics of the mosquito population, with rapid decay or growth. In fact, it seems that they are periods where the mosquito dynamics are mainly driven by temperature, while in other periods the dynamic is more driven by rainfall. Other approaches, based either on constant (model 4), temperature-dependent (model 2) or rainfall-dependent parameter values (model 3), provide similar results on average. However, they are not always satisfactory enough to be an "accurate" alternative, except model 3 when residual fertility is not too large.

Quality control within SIT is an important issue: if it fails, then release cannot occur. Within the quality control, we focus on residual fertility. We show that it may have an important impact on SIT duration and, eventually on its efficacy. In nuisance reduction, the lower the residual fertility, the lower the number of massive releases needed to switch from a massive releases strategy to a small releases strategy. However our results show that the massive releases duration for "high" residual fertility can greatly vary depending if the model takes into account the temperature, the rainfall or both. \yves{While the constraint $\varepsilon<1/\mathcal{N}$ is needed for elimination, it is not necessary to reduce the epidemiological risk.}

When SIT is considered, coupled or not with Mechanical control, we recommend starting the releases of sterile males within a period that last from July to December, when the mosquito population is, in general, at its lowest. However, the duration of the SIT treatment or the number of periodic releases, constant in size, of sterile males may vary thanks to the environmental parameters.

Clearly, combining Mechanical control with SIT is strongly recommended, in particular when the SIT treatment starts within a humid year: see for instance the years 2010 and 2011 where the dry period was more rainy than usual.

The massive and small releases strategy is only useful to reduce substantially the nuisance due to mosquitoes, i.e. to reach nearly elimination and to maintain the population under a certain level, related, here, with the number of sterile male individuals considered in the small releases \yves{strategy}. If the objective of SIT control is only to reduce the epidemiological risk, i.e. $\mathcal{R}_{eff}<0.5$, then the duration of the SIT treatment, with massive releases, will be short compared to elimination. In addition, having non-zero residual fertility, with $\varepsilon>1/\mathcal{N}$, seems to be less problematic to reach $\mathcal{R}_{eff}<0.5$, than to reach elimination. From a practical point of view, this can be very convenient. However, {once $\mathcal{R}_{eff}<0.5$ or elimination is reached,} it will be necessary to continue to release a sufficient amount of sterile males in order to keep the wild population under a certain threshold. At this stage, a closed-loop control (taking into account feedback from the system, like the size of the mosquito population through Mark Release and Recapture experiments) can be used from time to time, in order to reduce the overall cost, as described in \cite{Bliman2019,Aronna2020}. \yves{Whatever the objective, i.e. nuisance reduction or lower the epidemiological risk, it is better to have a residual fertility as small as possible in order to gain in the amount of sterile males to release, in the number of releases, and, thus, in sterile insects production.}

Our model, while minimalistic from the variables and parameters point of view, captures relatively well the dynamics of the \textit{Aedes albopictus} population throughout the year, without too many details. However, as with all models, improvements can be made. For instance, by taking into account the migration of males and females from neighboring places. Indeed, we have considered our area closed, in the sense that no external mosquitoes invade the treated area. Migration is another major concern in SIT treatment. We recently showed that if migration is small, SIT can be performed \cite{Bliman22}. Otherwise, it is necessary to isolate the target domain. This is not easy with mosquitoes because, so far, no (killing) attractant is efficient.

\yves{The main objective of SIT is to reduce the wild population, and in particular the wild females, because they are those who transmit diseases, like Dengue.} So far, in our model, we assumed that no sterile females are released. In fact, we know that this is not the case: there is always a percentage of sterile females that are released. \yves{In \cite{DUMONT2022}, outside an epidemic period, we show that releasing sterile females is not an issue. Within an epidemic period, releasing sterile females sterile is not an issue as long as the sterile females release rate is lower than a given threshold, $\Lambda_F^{crit}$, as estimated in \cite{DUMONT2022}. In particular, this sterile females threshold was estimated with constant parameters values and we know that it depends on $\mu_{A,2}$ and $\mathcal{R}_0^2$. Thus, according to our parameters variations, it is more than certain that the sterile females release rate threshold will vary greatly, such that the releases of sterile females will be more problematic in periods, when the wild population is small and the epidemiological risk is (very) low.} If, in the past, only less than $5\%$ of sterile females was acceptable \cite{Sharma1972}, this is not the case now. IAEA recommends not to release more than $1\%$ of sterile females. This is part of the control quality process. Thus, following \cite{DUMONT2022}, accidental releases of females could be taken into account, in order to derive how they could impact the SIT releases strategy, \yves{to reduce the epidemiological risk. Most certainly, there will be an impact in the duration of the releases, i.e. an increase, such that it will be necessary to adapt the releases strategy.}

\yves{Last, of course, it would be more than interesting to compare our SIT duration estimates with real (entomological and epidemiological) data recorded during a whole SIT campaign.}

\section*{Acknowledgments}

\yves{First of all, we would like to thank the anonymous reviewers for their remarks and advice that helped us improve the preliminary version of this paper.}

\noindent Preliminary results of this work have been presented during the workshop on SIT modeling, organized on R\'eunion island from the 27th of November to the 5th of December, 2021, with the support  of the TIS 2B project and the European Agricultural Fund for Rural Development (EAFRD). \yves{This work has been realized with the support of MESO@LR-Platform at the University of Montpellier.}
This work is partially supported by the "SIT feasibility project against \textit{Aedes albopictus} in R\'eunion Island", TIS 2B (2020-2022), jointly funded by the French Ministry of Health and the European Regional Development Fund (ERDF). YD is (partially) supported by the DST/NRF SARChI Chair in Mathematical Models and Methods in Biosciences and Bioengineering at the University of Pretoria (Grant 82770). YD acknowledges the support of DST/NRF Incentive Grant (Grant 119898). YD acknowledges the support of the Conseil R\'egional de la R\'eunion, the Conseil D\'epartemental de la R\'eunion, the European Agricultural Fund for Rural Development (EAFRD), the European Regional Development Fund (ERDF), and the Centre de Coop\'eration Internationale en Recherche Agronomique  pour le D\'eveloppement (CIRAD).

\noindent This paper is dedicated to the memory of Professor Svetoslav Markov (Bulgarian Academy of Sciences)
\bibliographystyle{abbrv}
\bibliography{bibliography.bib}

\section{Appendix A: Parameters values}
\label{AppendixA}
\label{tableparam}
See Table \ref{tab:parameters1}, page \pageref{tab:parameters1}, and Table \ref{tab:parameters2}, page \pageref{tab:parameters2}.

\begin{table}[]\small
\begin{center}
\caption{
{\bf Entomological parameters of \textit{Aedes albopictus} at different temperatures (from \cite{Delatte2009})}
}    
\label{tab:parameters1}
\begin{tabular}{|c|c|c|c|c|c|c|}
\hline 
Symbol & Name & T=15$^o$ & T=20$^o$ & T=25$^o$ & T=30$^o$ & T=35$^o$\tabularnewline
\hline 
\hline 
$r_{viable}$ & Proportion of viable eggs (E-L1) & $8.2$ & $66.9$ & $49.2$ & $51.4$ & $10$\tabularnewline
\hline 
$N_{eggs}$ & Number of eggs deposited & $0$ & $50.8$ & $65.3$ & $74.2$ & $48.7$\tabularnewline
\hline 
$\tau_{gono}$ & Duration of the gonotrophic cycle & NA & $8.1$ & $3.1$ & $3.9$ & $1.3$\tabularnewline
\hline 
$\tau_{A}$ & Time from hatching to emergence & $35$ & $14.4$ & $10.4$ & $8.8$ & $12.3$\tabularnewline
\hline 
$s_{A}$ & Survivorship from larva first instar to adult & $50$ & $77.5$ & $76.3$ & $67.5$ & $2.5$\tabularnewline
\hline 
$\tau_{M}$ & Adult male half-life  & $15.45$ & $10.25$ & $9.6$ & $8.55$ & $7.4$\tabularnewline
\hline 
$\tau_{F}$ & Adult female half-life & $19.65$ & $15.15$ & $15.3$ & $16.9$ & $10$\tabularnewline
\hline 
\end{tabular}
\end{center}
\end{table}
\begin{table}[]
   \caption{\bf Parameter values for system \eqref{eq:no SIT} deduced from Table \ref{tab:parameters1}, page \pageref{tab:parameters1}}
    \label{tab:parameters2}
{\small
\begin{tabular}{|c|c|c|c|c|c|c|c|}
\hline 
&&&mean&mean&mean&mean&mean\tabularnewline
Symbol&Name&Formula&value at&value at&value at&value at&value at\tabularnewline
&&&15$^o$&20$^o$&25$^o$&30$^o$&35$^o$\tabularnewline
\hline 
\hline 
$\phi$ & Effective fecundity & $\dfrac{r_{viable}N_{eggs}}{\tau_{gono}}$ & 0 & $4.1957$ & $10.3637$ & $9.7792$ & $3.7462$\tabularnewline
\hline 
$\mu_{A,1}$ & Aquatic death rate & $-\dfrac{log(s_{A})}{\tau_{A}}$ & $0.0198$ & $0.0177$ & $0.0260$ & $0.0447$ & $0.2999$\tabularnewline
\hline 
$\nu_{A}$ & Aquatic to adult  & $\dfrac{1}{\tau_{A}}$ & $0.0286$ & $0.0694$ & $0.0962$ & $0.1136$ & $0.0813$\tabularnewline
 & transition rate &  &  &  &  &  & \tabularnewline
\hline 
$\mu_{M}$ & Adult male death rate & $\dfrac{\log(2)}{\tau_{M}}$ & $0.0449$ & $0.0676$ & $0.0722$ & $0.0811$ & $0.0937$\tabularnewline
\hline 
$\mu_{F}$ & Adult female death rate & $\dfrac{log(2)}{\tau_{F}}$ & $0.0353$ & $0.0458$ & $0.0453$ & $0.0413$ & $0.0693$\tabularnewline
\hline 
\end{tabular}

\rewieverTwo{The formula given in Table \ref{tab:parameters2} are similar to those given in  \cite{Strugarek2019}. Let us explain how we deduce $\mu_M$ from $\tau_M$. Since $\tau_M$ is the adult-male half-life, it means that half of the population has disappeared in $\tau_M$ days. Then,  assuming that the dynamic of adult males follows an exponential law, i.e. $x(t)=x(0)e^{-\mu_M t}$, it is straightforward to deduce that $\mu_M$ is solution of $\dfrac{1}{2}=e^{-\mu_M\times \tau_M}$, that is $\mu_M=\dfrac{\log{2}}{\tau_M}$. The other parameter values follow the same reasoning.}
}
\end{table}



\section{Appendix B: Equilibria}
\label{AppendixB}
\yves{We are looking for the positive equilibria with and without SIT releases:}
\begin{itemize}
\item Without SIT releases.

In the case of constant coefficients, the equilibria of System \eqref{eq:no SIT} are $0$ and
\begin{equation}
    \left\{\begin{array}{l}
    A^*=\dfrac{\gamma+\mu_{A,1}}{\mu_{A,2}}(\mathcal{N}-1),\\
    M^*=\dfrac{(1-r)\gamma A^*}{\mu_M}=Q(\mathcal{N}-1),\\
    F^*=\dfrac{r\gamma A^*}{\mu_F},
    \end{array}\right.
    \label{eq_init}
\end{equation}
where $Q=\dfrac{(1-r)\gamma\left(\gamma+\mu_{1,A}\right)}{\mu_{2,A}\mu_{M}}$.
\item With SIT releases.
\label{SIT-appendix}

In the case of constant coefficients, some calculations are needed to find the release threshold, $M_{T_1,\varepsilon}$ and also to derive the equilibria. We have to solve 
\[
\left\{ \begin{array}{l}
\phi F=\left(\gamma+\mu_{1,A}+\mu_{2,A}A\right)A,\\
(1-r)\gamma A=\mu_{M}M,\\
\dfrac{M+\varepsilon\beta M_{T}}{M+\beta M_{T}}r\gamma A=\mu_{F}F.
\end{array}\right.
\]
Using the first and the third equalities leads to
\[
\dfrac{M+\varepsilon\beta M_{T}}{M+\beta M_{T}}r\gamma=\dfrac{\mu_{F}}{\phi}\left(\gamma+\mu_{1,A}+\mu_{2,A}A\right),
\]
and using the fact that 
\[
M=\dfrac{(1-r)\gamma}{\mu_{M}}A,
\]
we have
\[
\dfrac{(1-r)\gamma A+\varepsilon\mu_{M}\beta M_{T}}{(1-r)\gamma A+\mu_{M}\beta M_{T}}r\gamma=\dfrac{\mu_{F}}{\phi}\left(\gamma+\mu_{1,A}+\mu_{2,A}A\right),
\]
that is
\[
r\gamma\phi\left((1-r)\gamma A+\varepsilon\mu_{M}\beta M_{T}\right)=\mu_{F}\left((1-r)\gamma A+\mu_{M}\beta M_{T}\right)\left(\gamma+\mu_{1,A}+\mu_{2,A}A\right),
\]
leading to second order polynomial
\begin{multline*}
\mu_{F}(1-r)\gamma\mu_{2,A}A^{2}+\left[\mu_{F}(1-r)\gamma\left(\gamma+\mu_{1,A}\right)+\mu_{F}\mu_{M}\beta M_{T}\mu_{2,A}-r\gamma\phi(1-r)\gamma\right]A\\+\mu_{F}\mu_{M}\beta M_{T}\left(\gamma+\mu_{1,A}\right)-r\gamma\phi\varepsilon\mu_{M}\beta M_{T}=0
\end{multline*}
\begin{multline*}
(1-r)\gamma\mu_{2,A}A^{2}+\left[(1-r)\gamma\left(\gamma+\mu_{1,A}\right)+\mu_{M}\beta M_{T}\mu_{2,A}-\dfrac{r\gamma\phi}{\mu_{F}}(1-r)\gamma\right]A\\+\mu_{M}\beta M_{T}\left(\gamma+\mu_{1,A}\right)-\dfrac{r\gamma\phi}{\mu_{F}}\varepsilon\mu_{M}\beta M_{T}=0
\end{multline*}
\begin{multline*}
(1-r)\gamma\mu_{2,A}A^{2}+\left[(1-r)\gamma\left(\gamma+\mu_{1,A}\right)\left(1-\mathcal{N}\right)+\mu_{M}\beta M_{T}\mu_{2,A}\right]A\\+\mu_{M}\beta M_{T}\left(\gamma+\mu_{1,A}\right)\left(1-\varepsilon\mathcal{N}\right)=0,
\end{multline*}
or equivalently
\begin{equation*}
\dfrac{(1-r)\gamma}{\mu_{M}}A^{2}-\left[Q\left(\mathcal{N}-1\right)-\beta M_{T}\right]A+\beta M_{T}\dfrac{\left(\gamma+\mu_{1,A}\right)}{\mu_{2,A}}\left(1-\varepsilon\mathcal{N}\right)=0,
\end{equation*}
that is, using \eqref{eq_init},
\begin{equation}\label{eq: A}
\dfrac{(1-r)\gamma}{\mu_{M}}A^{2}-\left[M^*-\beta M_{T}\right]A+\beta M_{T}\dfrac{\left(\gamma+\mu_{1,A}\right)}{\mu_{2,A}}\left(1-\varepsilon\mathcal{N}\right)=0.
\end{equation}
We compute the discriminant of the last equation
\[
\Delta(\varepsilon)=\left[M^*-\beta M_{T}\right]^{2}-4Q\beta M_{T}\left(1-\varepsilon\mathcal{N}\right).
\]
We will distinguish three cases:
\begin{itemize}
    \item Assume $\varepsilon > 1/\mathcal{N}$, then $\Delta(\varepsilon)>0$, such that there always exists one positive equilibrium
$$
A_{\varepsilon}^{*}=\mu_M\dfrac{M^*-\beta M_{T}+\sqrt{\Delta(\varepsilon)}}{2(1-r)\gamma}.
$$
After straightforward computations, we can show that whatever the size of the massive releases, the aquatic equilibrium is bounded from below by
$$
2\dfrac{\varepsilon\mathcal{N}-1}{\mathcal{N}-1} A^*.
$$
\yves{All other values, $M_{\varepsilon}^*$ and $F_{\varepsilon}^*$, follow.}
\item Assume $\varepsilon=1/\mathcal{N}$. Then $\Delta(\varepsilon)=0$ iff $\beta M_{T}=M^*$. In fact if $\beta M_{T}>M^*$, then $A_{\varepsilon}^{*}=0$.
\item Assume $\varepsilon<1/\mathcal{N}$. Setting $y=\beta M_{T}$, we derive 
\[\begin{array}{rcl}
\Delta(\varepsilon)&=&\left[Q\left(\mathcal{N}-1\right)-y\right]^{2}-4Q\left(1-\varepsilon\mathcal{N}\right)y\\
&=&y^{2}-2Q\left(\left(\mathcal{N}-1\right)+2\left(1-\varepsilon\mathcal{N}\right)\right)y+\left(Q\left(\mathcal{N}-1\right)\right)^{2}
\end{array}\]
The discriminant of the equation $\Delta(\varepsilon)=0$ following the variable $y$ is given by
\[\begin{array}{rcl}
\delta(\varepsilon)&=&\left(2Q\right)^{2}\left(\left(\mathcal{N}-1\right)+2\left(1-\varepsilon\mathcal{N}\right)\right)^{2}-\left(2Q\left(\mathcal{N}-1\right)\right)^{2}\\
&=&16\left(1-\varepsilon\mathcal{N}\right)\left(1-\varepsilon\right)Q^{2}\mathcal{N}.
\end{array}\]
Since $\varepsilon<1/\mathcal{N}$, then $\delta(\varepsilon)>0$,
and we obtain two roots
\[
\beta M_{T_1,\varepsilon}=Q\left(\mathcal{N}+1-2\varepsilon\mathcal{N}-2\sqrt{\left(1-\varepsilon\mathcal{N}\right)\left(1-\varepsilon\right)\mathcal{N}}\right)
\]
and
\[
\beta M_{T_2,\varepsilon}=Q\left(\mathcal{N}+1-2\varepsilon\mathcal{N}+2\sqrt{\left(1-\varepsilon\mathcal{N}\right)\left(1-\varepsilon\right)\mathcal{N}}\right).
\]
When $\varepsilon=0$, we recover the result obtained in \cite{AnguelovTIS2020}
\[
\beta M_{T_{1},0}=Q\left(\mathcal{\sqrt{N}}-1\right)^{2}.
\]
Assume $0<\varepsilon<\dfrac{1}{\mathcal{N}}$. Then, when  $0<M_{T}<M_{T_{1},\varepsilon}$,
we have $\Delta(\varepsilon)>0$, and thus two equilibria
\[
A_{1,\varepsilon}=\mu_{M}\dfrac{Q\left(\mathcal{N}-1\right)-\beta M_{T}-\sqrt{\left[Q\left(\mathcal{N}-1\right)-\beta M_{T}\right]^{2}-4Q\beta M_{T}\left(1-\varepsilon\mathcal{N}\right)}}{2(1-r)\gamma},
\]
and 
\[
A_{2,\varepsilon}=\mu_{M}\dfrac{Q\left(\mathcal{N}-1\right)-\beta M_{T}+\sqrt{\left[Q\left(\mathcal{N}-1\right)-\beta M_{T}\right]^{2}-4Q\beta M_{T}\left(1-\varepsilon\mathcal{N}\right)}}{2(1-r)\gamma}.
\]
Since 
$$
\begin{array}{rcl}
A_{1,\varepsilon}+A_{2,\varepsilon}&=&\frac{\mu_M}{2(1-r)\gamma}\left(Q\left(\mathcal{N}-1\right)-\beta M_{T}\right)\\
&>&\frac{\mu_M}{2(1-r)\gamma}\left(Q\left(\mathcal{N}-1\right)-\beta M_{T_{1},\varepsilon}\right)\\
&=&\frac{\mu_M}{2(1-r)\gamma}\left(Q\left(2\varepsilon\mathcal{N}+2\left(\sqrt{\left(1-\varepsilon\mathcal{N}\right)\left(1-\varepsilon\right)\mathcal{N}}-1\right)\right)\right).
\end{array}$$
We remark that
\[
\left(1-\varepsilon\mathcal{N}\right)\left(1-\varepsilon\right)\mathcal{N}-1=\mathcal{N}-1-\left(\mathcal{N}+1\right)\varepsilon+\mathcal{N}\varepsilon^{2}.
\]
The discriminant of the last polynomial in $\varepsilon$ is given by
\[
\Delta_{\varepsilon}=\left(\mathcal{N}+1\right)^{2}-4\mathcal{N}\left(\mathcal{N}-1\right)=-3\mathcal{N}^2+6\mathcal{N}+1
\]
In general, $\N$ is large, i.e. $\mathcal{N}>>1+\frac{2}{\sqrt{3}}$, such that $\Delta_{\varepsilon}<0$, such that the $\varepsilon$-polynomial is always positive. 
Then $A_{1,\varepsilon}+A_{2,\varepsilon}>0$.
Moreover, since the last term in \eqref{eq: A}, 
we have also $A_{1,\varepsilon}A_{2,\varepsilon}>0$.
Thus the two roots $A_{1,\varepsilon}$ and $A_{2,\varepsilon}$ are positive. Therefore the system has two equilibria $\mathbf{E}_{1,2}=\left(A_{1,2},M_{1,2},F_{1,2}\right)$ such that
${\bf 0}<<\mathbf{E}_{1}<<\mathbf{E}_{2}$.
\end{itemize}
{In practice, we are, in general, in the case $\varepsilon <1/\N$, see \eqref{eq:1/N} and $\N>1+\frac{2}{\sqrt{3}}$.}
\end{itemize}

\paragraph{Appendix C: The epidemiological Model}
\label{AppendixC}
We briefly recall the SIR-SEI model of dengue transmission, studied in ~\cite{DUMONT2022}, without taking into account the accidental releases of sterile females.
The evolution of the human population is given by the following SIR model, with Susceptible, Infected, and Recovered compartments:
\begin{equation}
    \left\{
\begin{array}{l}
\dot S_h = \mu_h N_h - B\beta_{mh}(T) F_I \dfrac{S_h}{N_h} -\mu_h S_h,\\
\dot I_h = B\beta_{mh}(T) F_I \dfrac{S_h}{N_h} -(\eta_h+\mu_h) I_h,\\
\dot R_h = \eta_h I_h - \mu_h R_h,
\end{array}
\right.
\label{eq600}
\end{equation}

Extending the model of evolution of the mosquito population described in the previous sections, we use here a SEI model for the wild female mosquitoes, with Susceptible, Exposed and Infected compartments, adapted from \eqref{TIS_mod}:
\begin{equation}
    \left\{ \begin{array}{l}
\dfrac{dA}{dt}=\phi(T)(F_S+F_E+F_I)-\left(\gamma(T)+\mu_{1,A}(T)+\mu_{A,2}(T,R)A\right)A,\\ \\
\dfrac{dM}{dt}=(1-r(T))\gamma(T)A-\mu_{M}(T)M,\\\\
\dfrac{dF_S}{dt}=r(T)\gamma(T)\dfrac{M+\varepsilon \beta M_S}{M+\beta M_S}A- B\beta_{mh}(T) F_S \dfrac{I_h}{N_h}-\mu_{F}(T)F, \\ \\
\dot F_E = B\beta_{mh}(T) F_S \dfrac{I_h}{N_h} - (\nu_m(T)+\mu_F(T)) F_E\\ \\
\dot F_I = \nu_m(T) F_E - \mu_F(T) F_I.
\end{array}\right.
\end{equation}
where $M_S$ is driven by \eqref{eq:MS}. All parameters are described in the Table below
\begin{table}[h!]\small
\begin{center} 
  \caption{\bf Description of the epidemiological and entomological parameters}
 \begin{tabular}{|l|l|l|}
    \hline
    Symbol & Description & Unit\\
    \hline \hline \hline
    $1/\mu_{h}$ & Average lifespan of human & Day \\
    \hline
    $1/\eta_h$ & Average viremic period & Day \\
    \hline
    $B$ & Daily average mosquito bites & Day$^{-1}$ \\
    \hline
    $\beta_{mh}$ & Transmission probability from infected mosquito & - \\
    \hline
    $\beta_{hm}$ & Transmission probability from infected human & - \\
    \hline
    $\phi$ & Number of eggs at each deposit per capita   & Day$^{-1}$ \\
    \hline
    $\gamma$ & Maturation rate from larvae to adult  & Day$^{-1}$ \\
    \hline
    $\mu_{A,1}$ & Density independent mortality rate of the aquatic stage  & Day$^{-1}$ \\
    \hline
    $\mu_{A,2}$ & Density dependent mortality rate of the aquatic stage  & Day$^{-1}$ Individuals $^{-1}$ \\
    \hline
    $r$ & Sex ratio & - \\
     \hline
     $1/\nu_m$ & Average extrinsic incubation period (EIP) & Day \\
     \hline
    $1/\mu_F$ & Average lifespan of fertilized and eggs-laying females & Day \\
         \hline
    $1/\mu_M$ & Average lifespan of males & Day \\
         \hline
  \end{tabular}
  \end{center}
\end{table}
\end{document}